\RequirePackage{ifpdf}
\ifpdf % We are running pdfTeX in pdf mode
\documentclass[pdftex]{sigma}
\else
\documentclass{sigma}
\fi

\begin{document}

\allowdisplaybreaks

\renewcommand{\thefootnote}{$\star$}

\renewcommand{\PaperNumber}{085}

\FirstPageHeading

\ShortArticleName{String Functions for Af\/f\/ine Lie Algebras Integrable Modules}

\ArticleName{String Functions \\
for Af\/f\/ine Lie Algebras Integrable Modules\footnote{This paper is a
contribution to the Special Issue on Kac--Moody Algebras and Applications. The
full collection is available at
\href{http://www.emis.de/journals/SIGMA/Kac-Moody_algebras.html}{http://www.emis.de/journals/SIGMA/Kac-Moody{\_}algebras.html}}}

\Author{Petr KULISH~$^\dag$ and Vladimir LYAKHOVSKY~$^\ddag$}

\AuthorNameForHeading{P. Kulish and V. Lyakhovsky}

\Address{$^\dag$~Sankt-Petersburg Department of Steklov Institute of Mathematics,\\
\hphantom{$^\dag$}~Fontanka 27, 191023, Sankt-Petersburg, Russia}
\EmailD{\href{mailto:kulish@euclid.pdmi.ras.ru}{kulish@euclid.pdmi.ras.ru}}

\Address{$^\ddag$~Department of Theoretical Physics, Sankt-Petersburg State University,\\
\hphantom{$^\ddag$}~1 Ulyanovskaya Str., Petergof, 198904, Sankt-Petersburg, Russia}
\EmailD{\href{mailto:Vladimir.Lyakhovsky@pobox.spbu.ru}{Vladimir.Lyakhovsky@pobox.spbu.ru}}

\ArticleDates{Received September 15, 2008, in f\/inal form December 04,
2008; Published online December 12, 2008}

\Abstract{The recursion relations of branching coef\/f\/icients $k_{\xi }^{\left( \mu
\right) }$ for a module $L_{\frak{g}\downarrow \frak{h}}^{\mu }$ reduced
to a Cartan subalgebra $\frak{h}$ are transformed in order to place the
recursion shifts $\gamma \in \Gamma _{\frak{a}\subset \frak{h}}$
into the fundamental Weyl chamber. The new ensembles $F\Psi $
(the ``folded fans'') of shifts were constructed and the corresponding
recursion properties for the weights belonging to
the fundamental Weyl chamber were
formulated. Being considered simultaneously for the set of string functions
(corresponding to the same congruence class $\Xi _{v}$ of modules) the system
of recursion relations constitute an equation $\mathbf{M}_{\left( u\right)
}^{ \Xi _{v} } \mathbf{m}_{\left( u\right) }^{ \mu
 }=\mathbf{\delta }_{\left( u\right) }^{\mu }$ where the operator $%
\mathbf{M}_{\left( u\right) }^{ \Xi _{v} } $ 
is an invertible matrix
whose elements are def\/ined by the coordinates and
multiplicities of the shift weights in the
folded fans $F\Psi $ and the components of the vector $\mathbf{m}%
_{\left( u\right) }^{ \mu  }$ are the string function
coef\/f\/icients for $L^{\mu }$\ enlisted up to
an arbitrary f\/ixed grade $u$. The examples are
presented where the string functions for modules of $\frak{g}=A_{2}^{\left(
1\right) }$ are explicitly constructed demonstrating that the set of folded
fans provides a compact and ef\/fective tool to study the integrable highest
weight modules.}


\Classification{17B10; 17B67}

\section{Introduction}

We consider integrable modules $L^{\mu }$ with the highest weight $\mu $ for
af\/f\/ine Lie algebra $\frak{g}$ and are especially interested in the
properties of the string functions related to $L^{\mu }$. String functions
and branching coef\/f\/icients of the af\/f\/ine Lie algebras arise in the
computation of the local state proba\-bi\-lities for solvable models on square
lattice~\cite{DJKMO}. Irreducible highest weight modules with dominant
integral weights appear also in application of the quantum inverse
scattering method~\cite{LD} where solvable spin chains are studied in the
framework of the AdS/CFT correspondence conjecture of the super-string
theory (see~\cite{KAZA,BE} and references therein).

There are dif\/ferent ways to deal with string functions. One can use the BGG
resolution~\cite{BGG} (for Kac--Moody algebras the algorithm is described in
\cite{Kac,Wak1}), the Schur function series \cite{FauKing}, the BRST
cohomology \cite{Hwang}, Kac--Peterson formulas \cite{Kac} or the
combinatorial methods applied in \cite{FeigJimbo}.

Here we want to develop a new description for string functions by applying
the recursive formulas for weight multiplicities and branching coef\/f\/icients
obtained in \cite{IKL-1}.

It was proved in \cite{Kac} that for simply laced or twisted af\/f\/ine Lie
algebra and integrable mo\-du\-le~$L^{\mu }$ with the highest weight $\mu $ of
level $1$ the string function is unique:
\begin{equation*}
\sigma \big( e^{-\delta }\big) :=\prod_{n=1}^{\infty }\frac{1}{%
(1-e^{-n\delta })^{\mathrm{mult}(n\delta )}}. % \label{kac-string}
\end{equation*}
so that the corresponding formal character $\mathrm{ch}\left( L^{\mu
}\right) $ can be easily written down provided the set $\max (\mu )$ of
maximal weights for $L^{\mu }$ is known:
\begin{equation}
\mathrm{ch}\left( L^{\mu }\right) =\sigma \big( e^{-\delta }\big)
\sum_{\alpha \in M}e^{\mu +\alpha -\big( \frac{\left| \alpha \right| ^{2}}{2%
}+\left( \mu |\alpha \right) \big) \delta }  \label{kac-character-l-1}
\end{equation}
with
\[
M:=\left\{
\begin{array}{c}
\sum\limits_{i=1}^{r}\mathbf{Z}\alpha _{i}^{\vee }\text{ }\text{for untwisted
algebras or }A_{2r}^{\left( 2\right) } \\
\sum\limits_{i=1}^{r}\mathbf{Z}\alpha _{i}\text{ }\text{for }A_{r}^{\left( u\geq
2\right) }\text{ and }A\neq A_{2r}^{\left( 2\right) }
\end{array}
\right\}
\] (see also Corollary 2.1.6 in \cite{Wak2}). Comparing this expression with the Weyl--Kac formula
\[
\mathrm{ch}\left( L^{\mu }\right) =\frac{1}{R}\sum_{w\in W}\epsilon
(w)e^{w\circ (\mu +\rho )-\rho },
\]
where the character can be treated
as generated by the denominator $\frac{1}{R}$
acting on the set of singular vectors $\Psi ^{\left( \mu \right)
}=\sum\limits_{w\in W}\epsilon (w)e^{w\circ (\mu +\rho )-\rho }$
of the module $L^{\mu }$ we see that in the relation
(\ref{kac-character-l-1})
both factors on the right hand side are simplif\/ied: singular weights are
substituted by the maximal ones and instead of the factor $\frac{1}{R}$
the string function $\sigma _{u}\big( e^{-\delta }\big) $ is applied.

In this paper we shall demonstrate that similar transformations can be
def\/ined when the level~$k\left( \mu \right) $
is arbitrary. To f\/ind these transformations we use the
recursion properties of branching coef\/f\/icients $k_{\xi }^{\left( \mu \right)
}$ for the reduced module $L_{\frak{g}\downarrow \frak{a}}^{\mu }$ where
the subalgebra$\frak{a}$ has the same rank
as $\frak{g}$: $r(\frak{a})=r(\frak{g})$.
These properties are formulated in \cite{IKL-1} in terms of relations
\begin{equation*}
k_{\xi }^{\left( \mu \right) }=\sum_{\gamma \in \Gamma _{\frak{a}\subset
\frak{g}}}s\left( \gamma \right) k_{\xi +\gamma }^{\left( \mu \right)
}+\sum_{w\in W}\epsilon \left( w\right) \delta _{\xi ,\pi _{\frak{a}}\circ
\left( w\circ\left( \mu +\rho \right) -\rho \right) },  %\label{recursion-rel}
\end{equation*}
where $\pi _{\frak{a}}$\ is the projection to the weight space of $\frak{a}$
and $\Gamma _{\frak{a}\subset \frak{g}}$ is the fan of the injection $\frak{a%
}\longrightarrow \frak{g}$, that is the set of vectors def\/ined by the
relation
\[
1-\prod_{\alpha \in \left( \pi _{\frak{a}}\circ \Delta ^{+}\right)
}\left( 1-e^{-\alpha }\right) ^{\mathrm{{mult}\left( \alpha \right) -{mult}}%
_{\frak{a}}\mathrm{\left( \alpha \right) }}=
\sum_{\gamma \in \Gamma _{\frak{a}\subset \frak{g}}}%
s\left( \gamma \right) e^{-\gamma }
\]
(with $s\left( \gamma \right)\neq 0$).
In particular when $\frak{a}$ is a Cartan subalgebra $\frak{h}$ of $\frak{g}$
the coef\/f\/icients $k_{\xi }^{\left( \mu \right) }$ are just the multiplicities
of the weights of $L^{\mu }$ and the corresponding fan $\Gamma _{\frak{h}%
\subset \frak{g}}$ coincides with $\widehat{\Psi ^{\left( 0 \right) }}$~--
the set of singular weights $\psi \in P$ for the module $L^{0 }$.

In Section~\ref{section3} we demonstrate that this set can be ``folded'' $\widehat{\Psi
^{\left( 0 \right) }}\longrightarrow F\Psi $\ so that the new
shifts (the vectors of the folded fan) $f\psi \in F\Psi $
connect only the weights in the
closure of the fundamental Weyl chamber while the recursive property survives
in a new form. Thus the recursive relations are obtained for the coef\/f\/icients
of the string functions for the modules~$L^{\xi _{j}}$ whose highest weights
$\xi _{j}$ belong to the same congruence class $\Xi _{k;v}$. When these
relations are applied simultaneously to the set of string functions located
in the main Weyl chamber (Section~\ref{section4}) this results in the system of linear
equations for the string function coef\/f\/icients (collected in the vectors
$\mathbf{m}_{\left( s, u \right) }^{\left( \mu \right)}$).
This system can be written in a compact form $\mathbf{M}_{\left( u\right) }^{
\Xi _{v} } \mathbf{m}_{\left( u\right) }^{ \mu  }=
\mathbf{\delta }_{\left( u\right) }^{\mu }$ where the ope\-ra\-tor~$\mathbf{M}%
_{\left( u\right) }^{\Xi _{v} } $  is a matrix whose elements
are composed by the multiplicities of weights in the folded fans $%
F\Psi $. The set is solvable and the solution~-- the vector $%
\mathbf{m}_{\left( u\right) }^{ \mu  }$~-- def\/ines the string
functions for $L^{\mu }$\ up to an arbitrary minimal grade $u$. In the Section~\ref{section5}
some examples are presented where the string functions for modules of $%
\frak{g}=A_{2}^{\left( 1\right) }$ are explicitly constructed.

The set of folded fans provides a compact and ef\/fective method to construct
the string functions.

\section[Basic definitions and relations]{Basic def\/initions and relations}\label{section2}

Consider the af\/f\/ine Lie algebra $\frak{g}$ with the underlying
f\/inite-dimensional subalgebra $\overset{\circ }{\frak{g}}$.

The following notation will be used:

$L^{\mu }$  -- the integrable module of $\frak{g}$ with the highest weight
$\mu $;

$r$ -- the rank of the algebra $\frak{g}$;

$\Delta $ -- the root system;

$\Delta ^{+}$ -- the positive root system for $%
\frak{g}$;

$\mathrm{mult}\left( \alpha \right) $ -- the multiplicity of the root $%
\alpha $ in $\Delta $;

$\overset{\circ }{\Delta }$ -- the f\/inite root system of the subalgebra $%
\overset{\circ }{\frak{g}}$;

$\mathcal{N}^{\mu }$ -- the weight diagram of $L^{\mu }$;

$W$ -- the corresponding Weyl group;

$C^{\left( 0\right) }$ -- the fundamental Weyl chamber;

$\rho $ -- the Weyl vector;

$\epsilon \left( w\right) :=\det \left( w\right) $, $w \in W$;

$\alpha _{i}$ -- the $i$-th simple root for $\frak{g}$, $i=0,\ldots ,r$;

$\delta $ -- the imaginary root of $\frak{g}$;

$\alpha _{i}^{\vee }$ -- the simple coroot for $\frak{g}$, $i=0,\ldots ,r$;

$\overset{\circ }{\xi }$ -- the f\/inite (classical) part of the weight $\xi
\in P$;

$\lambda =\big( \overset{\circ }{\lambda };k;n\big) $ -- the
decomposition of an af\/f\/ine weight indicating the f\/inite part $\overset{\circ
}{\lambda }$, level $k$ and grade $n$;

$\overline{C_{k}^{\left( 0\right) }}$ -- the intersection of the closure of
the fundamental Weyl chamber $C^{\left( 0\right) }$\ with the plane with
f\/ixed level $k=\mathrm{const}$;

$P$ -- the weight lattice;

$Q$ -- the root lattice;

$M:=\left\{
\begin{array}{c}
\sum\limits_{i=1}^{r}\mathbf{Z}\alpha _{i}^{\vee }\text{ for untwisted algebras or }%
A_{2r}^{\left( 2\right) }, \\
\sum\limits_{i=1}^{r}\mathbf{Z}\alpha _{i}\text{ for }A_{r}^{\left( u\geq 2\right) }%
\text{ and }A\neq A_{2r}^{\left( 2\right) },
\end{array}
\right\} $;

$\mathcal{E}$  -- the group algebra of the group $P$;

$\Theta _{\lambda }:=e^{-\frac{\left| \lambda \right| ^{2}}{2k}\delta
}\sum\limits_{\alpha \in M}e^{t_{\alpha }\circ \lambda }$ -- the classical
theta-function;

$A_{\lambda }:=\sum\limits_{s\in \overset{\circ }{W}}\epsilon (s)\Theta
_{s\circ \lambda }$;

$\Psi ^{\left( \mu \right) }:=e^{\frac{\left| \mu +\rho \right| ^{2}}{2k}%
\delta \ -\ \rho }A_{\mu +\rho }=e^{\frac{\left| \mu +\rho \right| ^{2}}{2k}%
\delta \ -\ \rho }\sum\limits_{s\in \overset{\circ }{W}}\epsilon (s)\Theta
_{s\circ \left( \mu +\rho \right) }=$

\qquad $=\sum\limits_{w\in W}\epsilon (w)e^{w\circ(\mu +\rho )-\rho }$ --
the singular weight element for the $\frak{g}$-module $L^{\mu }$;

$\widehat{\Psi ^{\left( \mu \right) }}$ -- the set of singular weights $\psi
\in P$ for the module $L^{\mu }$ with the coordinates \newline
{}\hspace*{10mm} $\big( \overset{\circ
}{\psi },k,n,\epsilon \left( w\left( \psi \right) \right) \big) \mid
_{\psi =w\left( \psi \right) \circ (\mu +\rho )-\rho }$ (this set is similar
to $P_{\mathrm{nice}}^{\prime }\left( \mu \right) $ in \cite{Wak1});

$m_{\xi }^{\left( \mu \right) }$ -- the multiplicity of the weight $\xi \in
P $ in the module $L^{\mu }$;

$\mathrm{ch}\left( L^{\mu }\right) $ -- the formal character of $L^{\mu }$;

$\mathrm{ch}\left( L^{\mu }\right) =\frac{\sum\limits_{w\in W}\epsilon (w)e^{w\circ
(\mu +\rho )-\rho }}{\prod\limits_{\alpha \in \Delta ^{+}}\left( 1-e^{-\alpha
}\right) ^{\mathrm{{mult}\left( \alpha \right) }}}=\frac{\Psi ^{\left( \mu
\right) }}{\Psi ^{\left( 0\right) }}$ -- the Weyl--Kac formula;

$R:=\prod\limits_{\alpha \in \Delta ^{+}}\left( 1-e^{-\alpha }\right) ^{\mathrm{{%
mult}\left( \alpha \right) }}=\Psi ^{\left( 0\right) }$ -- the denominator;

$\max (\mu )$ -- the set of maximal weights of $L^{\mu }$;

$\sigma _{\xi }^{\mu }\left( q\right) =\sum\limits_{n=0}^{\infty }m_{\left( \xi
-n\delta \right) }^{\left( \mu \right) }q^{n}$ -- the string function
through the maximal weight $\xi $.

\section{Folding a fan}\label{section3}

The generalized Racah formula for
weight multiplicities  $m_{\xi
}^{\left( \mu \right) }$ (with $\xi \in P$)
in integrable highest weight modules $L^{\mu
}\left( \frak{g}\right) $ (see \cite{Ful} for a f\/inite dimensional variant),
\begin{equation}
m_{\xi }^{\left( \mu \right) }=-\sum_{w\in W\setminus e}\epsilon (w)m_{\xi
-\left( w\circ \rho -\rho \right) }^{\left( \mu \right) }+\sum_{w\in W}\epsilon
(w)\delta _{\left( w\circ(\mu +\rho )-\rho \right) ,\xi },  \label{gen-mult-form}
\end{equation}
can be obtained as a special case of developed in \cite{IKL-1} (see also~\cite{LDu}) branching algorithm
for af\/f\/ine Lie algebras. To apply this formula
(\ref{gen-mult-form}) we must determine two sets
of singular weights: $\widehat{\Psi ^{\left( \mu \right) }}$ for the module~$%
L^{\mu }$ and $\widehat{\Psi ^{\left( 0\right) }}$ for $L^{0}$. (As it was indicated
in the Introduction the set $\widehat{\Psi ^{\left( 0\right) }}$
coincides with the fan $\Gamma _{\frak{h}\subset \frak{g}}$  of the injection
$\frak{h}\longrightarrow \frak{g}$ of the Cartan subalgebra $\frak{h}$ in the Lie
algebra~$\frak{g}$.)

Our main idea is to contract the set $\widehat{\Psi ^{\left( 0\right) }}$
(the fan $\Gamma _{\frak{h}\subset \frak{g}}$) into the closure $\overline{%
C^{\left( 0\right) }}$ of the fundamental Weyl chamber $C^{\left( 0\right) }$.
We shall use the set $\max (\mu )$ of maximal weights of~$L^{\mu }\left(
\frak{g}\right) $ instead of~$\widehat{\Psi ^{\left( \mu \right) }}$. And as
a result we shall f\/ind the possibility to solve the relations based on the
recurrence properties of weight multiplicities, to obtain the explicit
expressions for the string functions $\sigma _{\xi \in \max (\mu )}^{\mu }$
and thus to describe the module~$L^{\mu }$.

Consider the module $L^{\mu}\left( \frak{g}\right) $ of level $k$: $\mu
=\big( \overset{\circ }{\mu };k;0\big) $. Let $\overline{C_{k;0}^{\left(
0\right) }}$ be the intersection of $\overline{C_{k}^{\left( 0\right) }}$
with the plane $\delta =0$, that is the ``classical'' part of the closure
of the af\/f\/ine Weyl chamber at level $k$.

To each $\xi \in P $ attribute a representative $w_{\xi }\in W$ of the class
of transformations
\[
w_{\xi }\in W/W_{\xi },\qquad W_{\xi }:=\left\{ w\in W\,|\,w\circ \xi =\xi
\right\} ,
\]
bringing the weight $\xi $ into the chamber $\overline{C_{k}^{\left(
0\right) }}$
\[
\bigg\{ w_{\xi }\circ \xi \in \overline{C_{k}^{\left( 0\right) }}\mid\xi \in
P,w_{\xi }\in W/W_{\xi }\bigg\} .
\]

Fix such representatives for each shifted vector
$\phi \left( \xi ,w\right) =\xi -\left( w\circ \rho
-\rho \right) $. The set
\[
\bigg\{ w_{\phi \left( \xi ,w\right) } \mid w_{\phi \left( \xi ,w\right)
}\circ \phi \left( \xi ,w\right) \in \overline{C_{k}^{\left( 0\right) }}
\bigg\},
\]
is in one-to-one correspondence with the set $\left\{ \phi \left( \xi
,w\right) \right\} $ of shifted weights. The recursion relation~(\ref{gen-mult-form}) can be written as
\begin{eqnarray*}
m_{\xi }^{\left( \mu \right) } &=&-\sum_{w\in W\setminus e}\epsilon
(w)m_{\phi \left( \xi ,w\right) }^{\left( \mu \right) }+\sum_{w\in
W}\epsilon (w)\delta _{\left( w\circ(\mu +\rho )-\rho \right) ,\xi } \\
&=&-\sum_{w\in W\setminus e}\epsilon (w)m_{w_{\phi \left( \xi ,w\right)
}\circ \phi \left( \xi ,w\right) }^{\left( \mu \right) }+\sum_{w\in
W}\epsilon (w)\delta _{\left( w\circ (\mu +\rho )-\rho \right) ,\xi }.
\end{eqnarray*}

Consider the restriction to $\overline{C_{k}^{\left( 0\right) }}$:
\begin{equation}
m_{\xi }^{\left( \mu \right) }\Big| _{\xi \in \overline{C_{k}^{\left(
0\right) }}}=-\sum_{w\in W\setminus e}\epsilon (w)m_{w_{\phi \left( \xi
,w\right) }\circ \phi \left( \xi ,w\right) }^{\left( \mu \right) }+\delta
_{\mu, \xi }.  \label{intermed-rel}
\end{equation}
In the r.h.s. the function $m_{\xi ^{\prime }}^{\left( \mu \right) }$ has an
argument $\xi ^{\prime }=w_{\phi \left( \xi ,w\right) }\circ \phi \left( \xi
,w\right) \in \overline{C_{k}^{\left( 0\right) }}$:
\[
m_{\xi ^{\prime }}^{\left( \mu \right) }=m_{w_{\phi \left( \xi ,w\right)
}\circ \phi \left( \xi ,w\right) }^{\left( \mu \right) }=m_{\xi +\left(
w_{\phi \left( \xi ,w\right) }\circ \phi \left( \xi ,w\right) -\xi \right)
}^{\left( \mu \right) }.
\]
Thus the new (``folded'') shifts are introduced:
\begin{gather*}
f\psi \left( \xi ,w\right)  := \left( \xi ^{\prime }-\xi \right) _{\xi
^{\prime }\neq \xi }=w_{\phi \left( \xi ,w\right) }\circ \left( \xi -\left(
w\circ \rho -\rho \right) \right) _{w\neq e}-\xi , \qquad
\xi ,\xi^{\prime }   \in   \overline{C_{k}^{\left( 0\right) }}, \qquad
\xi^{\prime } \neq \xi .
\end{gather*}
When the sum over $W\setminus e$ in the expression (\ref{intermed-rel}) is
performed the shifted weight $\xi ^{\prime }$ acquires the (f\/inite)
multiplicity $\widehat{\eta }\left( \xi ,\xi^{\prime }\right) $:
\begin{equation}
\widehat{\eta }\left( \xi ,\xi ^{\prime }\right) =-\sum_{ w\in
W\setminus e, } \epsilon (w),
\label{mult-in-ffan}
\end{equation}
(the sum is over all the elements $w\in W\setminus e $ satisfying the
relation $w_{\phi \left( \xi ,\tilde{w}\right) }\circ \left( \xi -\left(
w\circ \rho -\rho \right) \right) =\xi ^{\prime }$) such that
\begin{equation}
m_{\xi }^{\left( \mu \right) }\Big|_{\xi \in \overline{C_{k}^{\left( 0\right) }}%
}= \sum_{\xi ^{\prime }\in \overline{C_{k}^{\left( 0\right) }},\xi ^{\prime
}\neq \xi }\widehat{\eta }\left( \xi ,\xi ^{\prime }\right) m_{\xi +f\psi
\left( \xi ,w\right) }^{\left( \mu \right) }+\delta _{\xi ,\mu }.
\label{intermediate-2}
\end{equation}
The main property of the multiplicities $\widehat{\eta }\left( \xi ,\xi
^{\prime }\right) $
is that they do not depend directly on $n_{\xi}$.
%indirectly they depend on w_{\phi \left( \xi ,w\right) }

\begin{lemma}\label{lemma1}
Let $\psi =\rho -w\circ \rho $;  $\phi \left( \xi ,w\right) = \xi + \psi$;
  $\xi ^{\prime }:=w_{\phi \left( \xi
,w\right) }\circ \phi \left( \xi ,w\right)$; $\xi ,\xi ^{\prime }\in
\overline{C_{k}^{\left( 0\right) }}$. Then the corresponding folded shifts $%
f\psi \left( \xi ,w\right) = \xi^{\prime} -\xi$
and multiplicities $\widehat{\eta }\left( \xi
,\xi ^{\prime }\right) $ depend only on $k$, $\overset{\circ }{\xi }$, and~$w
$.
\end{lemma}

\begin{proof}
As far as imaginary roots are $W$-stable we have: $w_{\phi \left( \xi
,w\right) }\circ \left( \xi +\widetilde{n}\delta \right) =w_{\phi \left( \xi
,w\right) }\circ \xi +\widetilde{n}\delta $. Thus for both $\xi $ and $%
\widetilde{\xi }=\xi +\widetilde{n}\delta $ the representatives of the
classes bringing $\phi \left( \xi ,w\right) $ and $\phi \big( \widetilde{%
\xi },w\big) $ to the fundamental chamber $\overline{C_{k}^{\left(
0\right) }}$ can be taken equal: $w_{\phi \left( \xi ,w\right) }=w_{\phi
\left( \widetilde{\xi },w\right) }\;{\rm mod}W_{\xi }$. In the shift $f\psi
\left( \xi ,w\right) $ decompose the element $w_{\phi \left( \xi ,w\right)
}=t_{\phi \left( \xi ,w\right) }\cdot s_{\phi \left( \xi ,w\right) }$ into
the product of the classical ref\/lection $s_{\phi \left( \xi ,w\right) }$\ and the
translation $t_{\phi \left( \xi ,w\right) }$. Denote by $%
\theta _{\phi \left( \xi ,w\right) }^{\vee }$ the argument (belonging to $M$)
of the translation $t_{\phi \left( \xi ,w\right) }$. The direct computation
demonstrates that the weight $f\psi \left( \xi ,w\right) $ does not depend
on $n_{\xi }$:
\begin{gather*}
 f\psi \left( \xi ,w\right) =
 \left(
\begin{array}{c}
s_{\phi \left( \xi ,w\right) }\circ \big( \overset{\circ }{\xi }+\overset{%
\circ }{\psi }\big) -\overset{\circ }{\xi }+k\overset{\circ }{\theta }%
_{\phi \left( \xi ,w\right) }^{\vee },0, \vspace{1mm}\\
n_{-w\circ \rho }-\frac{k}{2}\big| \theta _{\phi \left( \xi ,w\right)
}^{\vee }\big| ^{2}-\big( s_{\phi \left( \xi ,w\right) }\circ \big(
\overset{\circ }{\xi }+\overset{\circ }{\psi }\big) ,\theta _{\phi \left(
\xi ,w\right) }^{\vee }\big)
\end{array}
\right) .  %\label{folded-shift}
\end{gather*}
Thus the shift $f\psi \left( \xi ,w\right) $ can be considered as depending
on $k$, $\overset{\circ }{\xi }$ and $w$: $f\psi =f\psi \big( \overset{%
\circ }{\xi },k,w\big) $. The multiplicity $\widehat{\eta }\left( \xi ,\xi
^{\prime }\right) $ (see (\ref{mult-in-ffan}))
depends only on the set of ref\/lections $w\in W$
connecting $\xi $ and $\xi ^{\prime }\neq \xi $
and does not depend on $n_{\xi }$ neither: $\widehat{\eta }\left( \xi ,\xi
^{\prime }\right) =\widehat{\eta }\big( \overset{\circ }{\xi },k,\xi
^{\prime }\big) $.
\end{proof}

Thus we have constructed the set of (nonzero) shifts $f\psi \big( \overset{%
\circ }{\xi },k,w\big) $ with the multiplicities $\widehat{%
\eta }\big( \overset{\circ }{\xi },k,\xi +f\psi \big( \overset{\circ }{\xi
},k,w\big) \big) $ and obtained the possibility to formulate the
recursion properties entirely def\/ined in the closure $\overline{%
C_{k}^{\left( 0\right) }}$ of the fundamental Weyl chamber.

Let us return to the relation (\ref{intermediate-2}),
\begin{gather*}
m_{\xi }^{\left( \mu \right) }\Big|_{\xi \in \overline{C_{k}^{\left( 0\right) }}%
}= \sum_{\xi ^{\prime }\in \overline{C_{k}^{\left( 0\right) }},\, \xi ^{\prime
}\neq \xi }\widehat{\eta }\big( \overset{\circ }{\xi },k,\xi ^{\prime
}\big) m_{\xi +f\psi  ( \overset{\circ }{\xi },k,w) }^{\left(
\mu \right) }+\delta _{\xi ,\mu }  \notag \\
\hphantom{m_{\xi }^{\left( \mu \right) }|_{\xi \in \overline{C_{k}^{\left( 0\right) }}}}{}
=\sum_{f\psi  ( \overset{\circ }{\xi },k,w ) \neq 0}
  \widehat{\eta }\big( \overset{\circ }{\xi }%
,k,\xi +f\psi \big( \overset{\circ }{\xi },k,w\big) \big) m_{\xi +f\psi
 ( \overset{\circ }{\xi },k,w ) }^{\left( \mu \right) }+\delta
_{\xi ,\mu }.  %\label{intermediate-3}
\end{gather*}
For simplicity from now on we shall omit some arguments and write down the
shifts as $f\psi \big( \overset{\circ }{\xi }\big) $ and their
multiplicities as $\widehat{\eta }\big( \overset{\circ }{\xi },\xi ^{\prime
}\big) $ (keeping in mind that we are at the level $k$ and the weight $\xi
^{\prime }$ depends on the initial ref\/lection $w$). The set of vectors:
\begin{gather*}
\widetilde{F\Psi }\big( \overset{\circ }{\xi }\big)
 := \left\{\xi^{\prime}- \xi =f\psi \big( \overset{\circ }{\xi }\big) =
\big( \overset{\circ }{f\psi }\big( \overset{\circ }{\xi }\big)
;0;n_{f\psi \big( \overset{\circ }{\xi }\big) }\big) \big|
\xi^{\prime} - \xi \neq 0  \right\} , \\
\xi ^{\prime }  =  w_{\phi \left( \xi ,w\right) }\circ \phi \left( \xi
,w\right),\qquad \xi ,\xi ^{\prime } \in \overline{C_{k}^{\left( 0\right) }},
\end{gather*}
plays here the role similar to that of the set $\left\{\Psi ^{\left(
0\right) }\setminus 0\right\}$ of nontrivial singular weights for $L^{0}$ in
the relation (\ref{gen-mult-form}) and is called the \emph{folded fan} for $%
\overset{\circ }{\xi }$. (The initial (unfolded) fan $\Gamma _{\frak{h}%
\subset \frak{g}}$ corresponds here to the injection of the Cartan
subalgebra.)

Thus we have proved the following property:

\begin{proposition}\label{proposition1}
Let $L^{\mu }$ be the integrable highest weight module of $\frak{g}$, $\mu
=\big( \overset{\circ }{\mu };k;0\big) $, $\xi =\big( \overset{\circ }{%
\xi };k;n_{\xi }\big) \in \mathcal{N}^{\mu }$, $\xi \in \overline{%
C_{k}^{\left( 0\right) }}$\ and let $\widetilde{F\Psi }\big( \overset{\circ
}{\xi }\big) $ be the folded fan for $\overset{\circ }{\xi }$ then
the multiplicity of the weight $\xi $ is subject to the recursion relation
\begin{equation}
m_{\xi }^{\left( \mu \right) }\Big|_{{\xi \in \overline{C_{k}^{\left(
0\right) }}}}=\sum_{f\psi  ( \overset{\circ }{\xi } ) \in
\widetilde{F\Psi } ( \overset{\circ }{\xi } ) }\widehat{\eta }%
\big( \overset{\circ }{\xi },\xi +f\psi \big( \overset{\circ }{\xi }%
\big) \big) m_{\xi +f\psi  ( \overset{\circ }{\xi } ) }^{\left(
\mu \right) }+\delta _{\xi ,\mu }.  \label{recursion-prop-weights}
\end{equation}
\end{proposition}

\section{Folded fans and string functions}\label{section4}

For the highest weight module $L^{\mu }\left( \frak{g}\right) $ with $\mu
=\big( \overset{\circ }{\mu };k;0\big) $ of level $k$ consider the set
of maximal vectors belonging to $\overline{C_{k}^{\left( 0\right) }}$
\[
{\cal Z} _{k}^{\mu}:=\bigg\{\zeta \in \max (\mu )\cap \overline{%
C_{k}^{\left( 0\right) }}  \bigg\}.
\]
Let $\pi $ be a projection to the subset of $P$ with level $k$ and grade $%
n=0$ and introduce the set:
\begin{equation*}
\Xi _{k}^{\mu }:=\big\{ \xi =\pi \circ \zeta  \mid \zeta \in
{\cal Z} _{k}^{\mu} \big\}.
%\label{maxes-for-mu}
\end{equation*}
The cardinality
\[
p_{\max }^{\left( \mu \right) }:=\#\left( \Xi _{k}^{\mu } \right)
\]
is f\/inite and we can enumerate the corresponding weights $\xi _{j}$:
\[
\Xi _{k}^{\mu } = \left\{ \xi _{j}\mid
j=1,\ldots ,p_{\max }^{\left( \mu \right) }\right\}.
\]
The string
functions necessary and suf\/f\/icient to construct the diagram $\mathcal{N}%
^{\mu }$\ (and correspondingly the character $\mathrm{ch}\left( L^{\mu
}\right) $) are
\begin{gather*}
\left\{ \sigma _{\zeta }^{\mu, k }\, |\, \zeta \in {\cal Z} _{k}^{\mu} \right\} ,
\qquad
\mathrm{ch}\left( L^{\mu }\right) =\sum_{\xi \in \max (\mu )}\sigma _{\xi
}^{\mu }\big( e^{-\delta }\big) e^{\xi }=\sum_{w \in W/ W_{\zeta}, \,
 \zeta \in {\cal Z} _{k}^{\mu}}
\sigma_{ \zeta }^{\mu, k }\big( e^{-\delta }\big) e^{w \circ \zeta }.
%\label{char-through-strings}
\end{gather*}
Let us consider these string functions as starting
from the points $\xi _{j}$ rather than from $\zeta$'s.
(For $\zeta = \xi_s -l\delta \in {\cal Z} _{k}^{\mu}$
the expansion
$\sigma _{\xi _{s}}^{\mu }\left( q\right)
=\sum\limits_{n=0}^{\infty }m_{\left( \xi _{s}-n\delta \right) }^{\left( \mu
\right) }q^{n}$ starts with string coef\/f\/icients
$m_{\left( \xi _{s}-n\delta \right)}|_{n<  l}=0$.)
Denote these extended string functions by
$\sigma _{j}^{\mu ,k}$ and introduce the set
\begin{equation*}
\Sigma _{k}^{\mu }:=\left\{ \sigma _{j}^{\mu ,k}\mid \xi _{j}\in \Xi
_{k}^{\mu }\right\} .  %\label{strings-for-mu}
\end{equation*}

Let us apply the relation (\ref{recursion-prop-weights}) to the weights of
the string $\sigma _{j}^{\mu ,k}\in \Sigma _{k}^{\mu }$ and
put $\xi =\xi _{j}+n_{j}\delta $,
\begin{gather*}
m_{ ( \overset{\circ }{\xi _{j}};k;n_{j} ) }^{\left( \mu \right) }=
 \sum_{f\psi  ( \overset{\circ }{\xi _{j}} ) \in \widetilde{F\Psi }%
 ( \overset{\circ }{\xi _{j}} ) }\widehat{\eta }\big( \overset{%
\circ }{\xi _{j}},\big( \overset{\circ }{\xi _{j}}+\overset{\circ }{f\psi }%
\big( \overset{\circ }{\xi_{j} }\big) ;k;n_{j}+n_{f\psi  ( \overset{%
\circ }{\xi _{j}} ) }\big) \big)  \\
\hphantom{m_{ ( \overset{\circ }{\xi _{j}};k;n_{j} ) }^{\left( \mu \right) }=}{}
\times m_{\big(  ( \overset{\circ }{\xi }_{j}+\overset{\circ }{f\psi }%
 ( \overset{\circ }{\xi_{j} } )  ) ;k; ( n_{j}+n_{f\psi
 ( \overset{\circ }{\xi_{j} } ) }) \big) }^{\left( \mu
\right) }+\delta _{\xi _{j},\mu }.
\end{gather*}
In the folded fan $\widetilde{F\Psi }\big( \overset{\circ }{\xi _{j}}%
\big) $ let us separate the summation over the grades $n_{f\psi }$ and the
classical parts $\overset{\circ }{f\psi }$ \ of the shifts $f\psi \big(
\overset{\circ }{\xi _{j}}\big) $. The overcrossing terms vanish because
their multiplicities are zero. The f\/irst term in the r.h.s.\
of the recursion relation takes the form
\begin{gather*}
  \sum_{n_{f\psi ( \overset{\circ }{\xi _{j}}) }}\sum
_{\substack{ \overset{\circ }{f\psi } ( \overset{\circ }{\xi_{j} }%
 ) ;  \\ f\psi  ( \overset{\circ }{\xi _{j}} ) \in \widetilde{%
F\Psi } ( \overset{\circ }{\xi _{j}} ) }}\widehat{\eta }\big(
\overset{\circ }{\xi _{j}},\big( \overset{\circ }{\xi _{j}}+\overset{\circ
}{f\psi }\big( \overset{\circ }{\xi_{j} }\big) ;k,;n_{j}+n_{f\psi  (
\overset{\circ }{\xi _{j}} ) }\big) \big)  m_{\big(  ( \overset{\circ }{\xi }_{j}+\overset{\circ }{%
f\psi } ( \overset{\circ }{\xi_{j} } )  ) ;k; (
n_{j}+n_{f\psi  ( \overset{\circ }{\xi_{j} } ) } ) \big)
}^{\left( \mu \right) }.
\end{gather*}
For the same reason we can spread the f\/irst summation over all the positive
grades. It is suf\/f\/icient to include the vector
with zero coordinates into the folded fan and put the multiplicity
$\eta \big( \overset{\circ }{\xi },\xi \big) =-1$. Introduce the set
\[
F\Psi \big( \overset{\circ }{\xi _{j}}\big) :=\widetilde{F\Psi }\big(
\overset{\circ }{\xi _{j}}\big) \cup \left( 0;0;0\right) .
\]
It is called \emph{the full folded fan}
or simply \emph{the folded fan} when
from the context it is clear what fan $\widetilde{F\Psi }\big( \overset{%
\circ }{\xi _{j}}\big) $ or $F\Psi \big( \overset{\circ }{\xi _{j}}%
\big) $ is actually used. The set of
multiplicities $\eta \big( \overset{\circ }{\xi },\xi ^{\prime }\big) $
for the shifts in $F\Psi \big( \overset{\circ }{\xi }\big) $ is thus
f\/ixed as follows:
\begin{equation}
\eta \big( \overset{\circ }{\xi },\xi ^{\prime }\big)\big|_{\xi ^{\prime
}-\xi \in F\Psi ( \overset{\circ }{\xi }) }:=-\sum_{\substack{ %
w\in W,  \\ w_{\phi ( \xi ,w) }\circ ( \xi -( w\circ
\rho -\rho ) ) =\xi ^{\prime }}}\epsilon (w),
\label{full-fan-multiplicity}
\end{equation}
and the recursion property (\ref{recursion-prop-weights}) is reformulated:
\begin{equation*}
\sum_{f\psi  ( \overset{\circ }{\xi } ) \in F\Psi  ( \overset{%
\circ }{\xi } ) }\eta \big( \overset{\circ }{\xi },\xi +f\psi \big(
\overset{\circ }{\xi }\big) \big) m_{\xi +f\psi  ( \overset{\circ }{%
\xi } ) }^{\left( \mu \right) }+\delta _{\xi ,\mu }=0, \qquad
\xi \in \overline{C_{k}^{\left( 0\right) }}.
%\label{recursion-prop-zero}
\end{equation*}
For the string $\sigma_{j}^{\mu,k }$
we can rewrite this relation separating the summations:
\begin{gather*}
 \sum_{n=0}^{\infty }\sum_{\substack{ \overset{\circ }{f\psi } (
\overset{\circ }{\xi_{j} } )  \\ f\psi  ( \overset{\circ }{\xi _{j}}%
 ) \in F\Psi  ( \overset{\circ }{\xi _{j}} ) }}\eta \big(
\overset{\circ }{\xi _{j}},\big( \overset{\circ }{\xi _{j}}+\overset{\circ
}{f\psi }\big( \overset{\circ }{\xi_{j} }\big) ;k;n_{j}+n\big) \big)
  m_{ (  ( \overset{\circ }{\xi }_{j}+\overset{\circ }{%
f\psi } ( \overset{\circ }{\xi_{j} } )  ) ;k; (
n_{j}+n )  ) }^{\left( \mu \right) }+\delta _{\xi _{j},\mu }=0.
\end{gather*}
The properties of $\mathcal{N}^{\mu }$ for an integrable
modules $L^{\mu }$ guarantee that for any f\/inite $n_{j}$ the f\/irst sum is
f\/inite. It extends to $n\leq -n_{j}$ (remember that $n_{j}$ is negative).
The second sum can also be augmented so that the vectors $ \big(
\overset{\circ }{\xi _{j}}+\overset{\circ }{f\psi }\big( \overset{\circ }{%
\xi_{j} }\big);k;0 \big) =\big( \overset{\circ }{\xi _{s}};k;0 \big)$
run over the set $\Xi_{k}^{\mu }$. Now
taking into account that $n_{j,s}$ does not depend on $n_j$ (Lemma~\ref{lemma1})
the notation can be simplif\/ied:
\begin{gather*}
\eta _{j,s}\left( n\right)  : =\eta \big( \overset{\circ }{\xi _{j}}%
,\big( \overset{\circ }{\xi _{s}};k;n_{j}+n\big) \big) , \qquad
m_{s,n_{j}+n}^{\left( \mu \right) }  : =m_{ ( \overset{\circ }{\xi _{s}}%
;k;n_{j}+n ) }^{\left( \mu \right) },
\end{gather*}
and the recursion property for the string functions in
$\big\{ \sigma _{j}^{\mu }|\xi
_{j}\in \Xi _{k}^{\mu }\big\} $ can be stated:

\begin{proposition}\label{proposition2}
Let $L^{\mu }$ be the integrable highest weight module of $\frak{g}$, $\mu
=\big( \overset{\circ }{\mu };k;0\big) $, $p_{\max }^{\left( \mu \right)
}:=\#\left( \Xi_{k}^{\mu }\right) $\ ,
$\xi _{j}=\big( \overset{\circ }{\xi _{j}};k;n_{j}\big) \in \Xi _{k}^{\mu }
+ n_j \delta $ , let $F\Psi \big( \overset{\circ }{\xi _{j}}\big) $ be the full folded
fan for $\overset{\circ }{\xi _{j}}$ and
$\eta _{j,s}\left( n\right)
=-\sum\limits_{ \tilde{w}_{j,s}, }\epsilon (\tilde{w}_{j,s})$
where the summation is over the elements $\tilde{w}_{j,s}$
of $W$ satisfying the equation
$ w_{\phi \left( \xi ,w\right) }\circ \left( \xi
_{j}-\left( \tilde{w}_{j,s}\circ \rho -\rho \right) \right) =\big( \overset{\circ }{\xi
_{s}};k;n_{j}+n\big),$
then for the string function coefficients
$m_{s,n_{j}+n}^{\left( \mu \right) }$
the following relation holds:
\begin{equation}
\sum_{s=1}^{p_{\max }^{\left( \mu \right) }}\sum_{n=-n_{j}}^{0}\eta
_{j,s}\left( n\right) m_{s,n_{j}+ n}^{\left( \mu \right) }=-\delta
_{\xi _{j},\mu }.  \label{recursion-prop-string}
\end{equation}
\end{proposition}

For a f\/ixed $n_{j}\leq 0$ consider the sequence of the string weights
\[
\xi _{j;n_{j}}=\big( \overset{\circ }{\xi _{j}};k;n_{j}\big) , \quad \xi
_{j;n_{j}+1}=\big( \overset{\circ }{\xi _{j}};k;n_{j}+1\big) ,\quad \ldots ,\quad \xi
_{j;0}=\big( \overset{\circ }{\xi _{j}};k;0\big) ,
\]
and write down two $\left( \left| n_{j}\right| +1\right) $-dimensional
vectors: the coordinates of the f\/irst one are the coef\/f\/icients of the $s$-th string $\left\{ \sigma _{s}^{\mu }\right\} $,
\begin{equation*}
\mathbf{m}_{\left( s;n_{j}\right) }^{\left( \mu \right) }:=\left(
m_{s,n_{j}}^{\left( \mu \right) },m_{s,n_{j}+1}^{\left( \mu \right) },\ldots
,m_{s,0}^{\left( \mu \right) }\right) ,
%\label{string-weights-multiplicities}
\end{equation*}
the second indicates that the $j$-th string $\sigma _{j}^{\mu ,k}$\ is
starting at the highest weight $\mu $,
\[
\mathbf{\delta }_{\left( j;n_{j}\right) }^{\mu }:=\left( 0,0,\ldots
,-1\right) .
\]
For the weights with $n\geq n_{j}$ we have the sequence of relations of the
type (\ref{recursion-prop-string}):
\begin{gather}
\sum_{s=1}^{p_{\max }^{\left( \mu \right) }}\sum_{n=0}^{-n_{j}}\eta
_{j,s}\left( n\right) m_{s,n_{j}+n}^{\left( \mu \right) }  = 0,  \notag
\\
\sum_{s=1}^{p_{\max }^{\left( \mu \right) }}\sum_{n=0}^{-n_{j}-1}\eta
_{j,s}\left( n\right) m_{s,n_{j}+n+1}^{\left( \mu \right) }  = 0,
\notag \\
  \cdots \cdots \cdots \cdots \cdots \cdots \cdots \cdots \cdots\cdots ,  \notag \\
\sum_{s=1}^{p_{\max }^{\left( \mu \right) }}\eta _{j,s}\left( 0\right)
m_{s,0}^{\left( \mu \right) }  = -1.  \label{preliminary-set-equations}
\end{gather}
Introduce the upper triangular $\left( \left| n_{j}\right| +1\right) \times
\left( \left| n_{j}\right| +1\right) $-matrix
\begin{equation*}
\mathbf{M}_{\left( j,s\right) }^{\Xi \mu}:=
\begin{array}{cccc}
\eta _{j,s}\left( 0\right)  & \eta _{j,s}\left( 1\right)  & \cdots
& \eta _{j,s}\left( -n_j \right)  \\
0 & \eta _{j,s}\left( 0 \right)  & \cdots  & \eta _{j,s}\left(-n_j -1\right)
\\
\vdots  & \vdots  & \vdots  & \vdots  \\
0 & 0 & \cdots  & \eta _{j,s}\left( 0\right)
\end{array}
.  %\label{blocks-of-etas}
\end{equation*}
The set of relations (\ref{preliminary-set-equations}) reads:
\begin{equation}
\mathbf{M}_{\left( j,s\right) }^{\Xi \mu} \cdot%
\mathbf{m}_{\left( s;n_{j}\right) }^{\left( \mu \right) }=\mathbf{\delta }%
_{\left( j;n_{j}\right) }^{\mu }.  \label{preliminary-matrix-relations}
\end{equation}
Perform the same procedure for the other weights $\xi _{j}\in \Xi _{k}^{\mu }$
putting the minimal values of grade equal: $n_{j}|_{j=1,\ldots ,p_{\max
}^{\left( \mu \right) }}=u$, that is construct all the folded fans $F\Psi
\big( \overset{\circ }{\xi _{j}}\big) $ (till the grade $u$) and the
corresponding sets of multiplicities $\eta _{j,s}\left( n\right) $
(def\/ined by relations (\ref{full-fan-multiplicity})). For $j=1,\ldots
,p_{\max }^{\left( \mu \right) }$ compose $\left( |u| +1\right)^2$ equations
of the type (\ref{preliminary-matrix-relations}):
\begin{gather}
\mathbf{M}_{\left( j,s\right) }^{\Xi \mu}
\mathbf{m}_{\left( s;n_j \right) }^{\left( \mu \right) }  = \mathbf{\delta }%
_{\left( j;n_j \right) }^{\mu },  \qquad
j,s  =  1,\dots, p_{\max }^{\left( \mu \right)}.
\label{almost-final-set-relations}
\end{gather}
Form two $\left( \left| u\right| +1\right) \times p_{\max }^{\left( \mu
\right) }$-dimensional vectors: the f\/irst with the string coef\/f\/icients,
\begin{gather*}
\mathbf{m}_{\left( u\right) }^{\left( \mu \right) }:=  \bigg(
m_{1,u}^{\left( \mu \right) },m_{1,u+1}^{\left( \mu \right) },\ldots
,m_{1,0}^{\left( \mu \right) },m_{2,u}^{\left( \mu \right)
},m_{2,u+1}^{\left( \mu \right) },\ldots ,m_{2,0}^{\left( \mu \right)
},\ldots    \notag \\
\hphantom{\mathbf{m}_{\left( u\right) }^{\left( \mu \right) }:=\bigg(}{}  \ldots ,  m_{p_{\max }^{\left( \mu \right) },u}^{\left( \mu \right)
},m_{p_{\max }^{\left( \mu \right) },u+1}^{\left( \mu \right) },\ldots
,m_{p_{\max }^{\left( \mu \right) },0}^{\left( \mu \right) }\bigg) ,
\end{gather*}
the second indicating that the string $\sigma _{j}^{\mu ,k}$ with number
$j$ starts at the highest weight $\mu $,
\[
\mathbf{\delta }_{\left( u\right) }^{\mu }:=\left( 0,0,\ldots ,0,0,0,\ldots
,0,0,0,\ldots ,-1,0,0,\ldots ,0\right) ,
\]
(here only in the j-th subsequence the last ($\left( \left| u\right|
+1\right) $-th) coordinate is not zero). Def\/ine the $\left( \left| u\right|
+1\right) p_{\max }^{\left( \mu \right) }\times \left( \left| u\right|
+1\right) p_{\max }^{\left( \mu \right) }$-matrix -- the block-matrix with
the blocks $\mathbf{M}_{\left( j,s\right) }^{\Xi \mu}$:
\[
\mathbf{M}^{\Xi \mu}:=\left\| \mathbf{M}_{\left( j,s\right) }^{\Xi \mu}\right\| _{j,s=1,\ldots,
p_{\max }^{\left( \mu \right) }}.
\]
In these terms the relations (\ref{almost-final-set-relations}) have the
following integral form:
\begin{equation}
\mathbf{M}^{\Xi \mu}\,\,\mathbf{m}_{\left( u\right) }^{\left( \mu \right) }=\mathbf{%
\delta }_{\left( u\right) }^{\mu }.  \label{final-equation}
\end{equation}
The matrix $\mathbf{M}^{\Xi \mu}$ being invertible
the equation (\ref{final-equation}) can be
solved. Thus we have demonstrated that the strings $\sigma _{j}^{\mu ,k}$
are determined by the matrix $\mathbf{M}^{\Xi \mu}$ whose elements are the full folded
fan weight multiplicities:

\begin{proposition}\label{proposition3}
Let $L^{\mu }$ be an integrable highest weight module of $\frak{g}$, $\mu
=\big( \overset{\circ }{\mu };k;0\big) $, $p_{\max }^{\left( \mu \right)
}:=\#\left( \Xi_{k}^{\mu}\right) $, $\xi _{j}=\big( \overset{\circ }{\xi _{j}};k;n_{j}\big) \in \Xi _{k}^{\mu }
+ n_j \delta
$; let $F\Psi \big( \overset{\circ }{\xi _{j}}\big) $ be the full folded
fan for $\overset{\circ }{\xi _{j}}$ and $\mathbf{M}^{\Xi \mu}$~-- the $\left( \left|
n_{j}\right| +1\right) p_{\max }^{\left( \mu \right) }\times \left( \left|
n_{j}\right| +1\right) p_{\max }^{\left( \mu \right) }$-matrix formed by the
blocks $\mathbf{M}_{\left( j,s\right) }^{\Xi \mu}$
\[
\mathbf{M}_{\left( j,s\right) }^{\Xi \mu}:=
\begin{array}{cccc}
\eta _{j,s}\left( 0 \right)  & \eta _{j,s}\left( 1\right)  & \cdots
& \eta _{j,s}\left( -n_{j} \right)  \\
0 & \eta _{j,s}\left( 0 \right)  & \cdots  & \eta _{j,s}\left( -n_{j} -1\right)
\\
\vdots  & \vdots  & \vdots  & \vdots  \\
0 & 0 & \cdots  & \eta _{j,s}\left( 0\right)
\end{array}
\]
where the elements $\eta _{j,s}\left( n\right) $ are the
multiplicities of the folded fan weights,
\[
\eta _{j,s}\left( n\right)
=-\sum_{ \tilde{w}_{j,s}, }\epsilon (\tilde{w}_{j,s})
\]
with the summation over the elements $\tilde{w}_{j,s} \in W$
satisfying the equation
\[
w_{\phi \left( \xi ,w\right) }\circ \left( \xi
_{j}-\left( \tilde{w}_{j,s}\circ \rho -\rho \right) \right) =\big( \overset{\circ }{\xi
_{s}};k;n_{j}+n\big).
\]
Let the string function coefficients
be the coordinates in the $n_{j}+1$-subsequences of the vector~$\mathbf{m}_{\left( n_{j}\right) }^{\left( \mu \right) }$.
Then for the coefficients of
$\big\{\sigma _{j}^{\mu ,k}\mid j=1,\ldots,
p_{\max }^{\left( \mu \right)}\big\}$ the following relation folds:
\begin{gather}
\mathbf{m}_{\left( n_{j}\right) }^{\left( \mu \right) }=
\left(\mathbf{M}^{\Xi \mu}\right)^{-1}
\mathbf{\delta }_{\left( n_{j}\right) }^{\mu }.
\label{finite-string-solution}
\end{gather}
\end{proposition}

Thus the solution $\mathbf{m}_{\left( n_{j} \right) }^{\left( \mu \right) }$
describes all the string functions relevant to the chosen module $L^{\mu }$
(with the grades no less than the preliminary f\/ixed $n_{j}=u$).

To describe the complete string functions it is suf\/f\/icient to send $u$ to
the limit $u\rightarrow -\infty $ .

\section{Examples}\label{section5}

\subsection[$\frak{g}=A_{2}^{\left( 1\right) }$]{$\boldsymbol{\frak{g}=A_{2}^{\left( 1\right) }}$}\label{section5.1}

Consider the fan $\Gamma _{\frak{h}\subset \frak{g}}$ (with $n_{\psi
^{\left( 0\right) }}\leq 9$):
\begin{gather}
\Gamma _{\frak{h}\subset \frak{g}}  = \big\{   \left(
0, 1, 0, 0, 1\right) , \left( 2, 1, 0, 0, -1\right) , \left(
1, 0, 0, 0, 1\right) , \left( 1, 2, 0, 0, -1\right) ,  \notag \\
\hphantom{\Gamma _{\frak{h}\subset \frak{g}}  = \big\{}{} \left( 2, 2, 0, 0, 1\right) ,\left( 3, 1, 0, 1, 1\right) , \left(
-1, 1, 0, 1, -1\right) ,\left( 1, 3, 0, 1, 1\right) ,   \notag \\
\hphantom{\Gamma _{\frak{h}\subset \frak{g}}  = \big\{}{}
\left( 1, -1, 0, 1, -1\right) , \left( 3, 3, 0, 1, -1\right)
,\left( -1, -1, 0, 1, 1\right) ,\left( 3, 4, 0, 2, 1\right) ,
\notag \\
\hphantom{\Gamma _{\frak{h}\subset \frak{g}}  = \big\{}{}
\left( 0, -2, 0, 2, 1\right) , \left( 2, 4, 0, 2, -1\right)
, \left( -1, -2, 0, 2, -1\right) ,\left( 4, 3, 0, 2, 1\right) ,
\notag \\
\hphantom{\Gamma _{\frak{h}\subset \frak{g}}  = \big\{}{}
\left( -2, 0, 0, 2, 1\right) , \left( 4, 2, 0, 2, -1\right)
, \left( -2, -1, 0, 2, -1\right) , \left( 0, 3, 0, 2, -1\right) ,
\notag \\
\hphantom{\Gamma _{\frak{h}\subset \frak{g}}  = \big\{}{}
\left( 3, 0, 0, 2, -1\right) , \left( -1, 2, 0, 2, 1\right)
, \left( 2, -1, 0, 2, 1\right) , \left( 0, 4, 0, 4, 1\right) ,
\notag \\
\hphantom{\Gamma _{\frak{h}\subset \frak{g}}  = \big\{}{}
\left( -3, -2, 0, 4, 1\right) ,\left( 5, 4, 0, 4, -1\right)
, \left( 2, -2, 0, 4, -1\right) , \left( 4, 0, 0, 4, 1\right) ,
\notag \\
\hphantom{\Gamma _{\frak{h}\subset \frak{g}}  = \big\{}{}
\left( -2, -3, 0, 4, 1\right) , \left( 4, 5, 0, 4, -1\right) ,
\left( -2, 2, 0, 4, -1\right) , \left( -3, 0, 0, 4, -1\right) ,
\notag \\
\hphantom{\Gamma _{\frak{h}\subset \frak{g}}  = \big\{}{}
\left( 1, -3, 0, 5, 1\right) , \left( 5, 1, 0, 5, -1\right)
, \left( 5, 5, 0, 5, 1\right) , \left( 1, 5, 0, 5, -1\right) ,
\notag \\
\hphantom{\Gamma _{\frak{h}\subset \frak{g}}  = \big\{}{}
\left( 0, -3, 0, 4, -1\right) , \left( 2, 5, 0, 4, 1\right)
, \left( 5, 2, 0, 4, 1\right) , \left( -3, -3, 0, 5, -1\right) ,
\notag \\
\hphantom{\Gamma _{\frak{h}\subset \frak{g}}  = \big\{}{}\left( -3, 1, 0, 5, 1\right) , \left( 6, 4, 0, 6, 1\right)
, \left( 3, -2, 0, 6, 1\right) , \left( -1, 4, 0, 6, -1\right) ,
\notag \\
\hphantom{\Gamma _{\frak{h}\subset \frak{g}}  = \big\{}{}
\left( -4, -2, 0, 6, -1\right) ,\left( 4, 6, 0, 6, 1\right)
, \left( -2, 3, 0, 6, 1\right) , \left( 4, -1, 0, 6, -1\right) ,
\notag \\
\hphantom{\Gamma _{\frak{h}\subset \frak{g}}  = \big\{}{}
\left( -2, -4, 0, 6, -1\right) , \left( 3, 6, 0, 6, -1\right)
,\left( 6, 3, 0, 6, -1\right) , \left( -4, -1, 0, 6, 1\right) ,
\notag \\
\hphantom{\Gamma _{\frak{h}\subset \frak{g}}  = \big\{}{}
\left( -1, -4, 0, 6, 1\right) , \left( 6, 1, 0, 8, 1\right)
, \left( -4, 1, 0, 8, -1\right) ,\left( 1, 6, 0, 8, 1\right) ,
\notag \\
\hphantom{\Gamma _{\frak{h}\subset \frak{g}}  = \big\{}{}
\left( 1, -4, 0, 8, -1\right) , \left( 6, 6, 0, 8, -1\right)
, \left( -4, -4, 0, 8, 1\right) , \left( 3, 7, 0, 9, 1\right) ,
\notag \\
\hphantom{\Gamma _{\frak{h}\subset \frak{g}}  = \big\{}{}
\left( -3, -5, 0, 9, 1\right) , \left( 5, 7, 0, 9, -1\right)
, \left( -1, -5, 0, 9, -1\right) , \left( 7, 3, 0, 9, 1\right) ,
\notag \\
\hphantom{\Gamma _{\frak{h}\subset \frak{g}}  = \big\{}{}
\left( -5, -3, 0, 9, 1\right) ,\left( 7, 5, 0, 9, -1\right)
, \left( -5, -1, 0, 9, -1\right) , \left( -3, 3, 0, 9, -1\right)
,   \notag \\
\hphantom{\Gamma _{\frak{h}\subset \frak{g}}  = \big\{}{}
\left( 3, -3, 0, 9, -1\right) , \left( -1, 5, 0, 9, 1\right) ,
\left( 5, -1, 0, 9, 1\right),   \ldots \big\} .  \label{star-a2-1}
\end{gather}
Here the f\/irst two coordinates are classical in the basis of simple roots $%
\left\{ \alpha _{1},\alpha _{2}\right\} $, next comes the level $k=0$, the
grade $n_{\psi ^{\left( 0\right) }}$ and the multiplicity $m_{\psi ^{\left(
0\right) }}$ of the weight $\psi ^{\left( 0\right) }\in \Gamma _{\frak{h}%
\subset \frak{g}}$ (for the injection $\frak{h} \longrightarrow \frak{g}$
we have $m_{\psi ^{\left( 0\right) }}=-\epsilon(w)$).

\subsubsection[$k=1$]{$\boldsymbol{k=1}$}\label{section5.1.1}

The set $\overline{C_{1;0}^{\left( 0\right) }}$ contains three weights ($%
p_{\max }^{\left( \mu \right) }=3$):
\begin{gather*}
\overline{C_{1;0}^{\left( 0\right) }}  = \big\{ \left( 0,0;1;0\right)
, ( \overset{\circ }{\omega _{1}};1;0 ) , ( \overset{\circ }{%
\omega _{2}};1;0 ) \big\} = \left\{\omega _{0},\omega _{1},\omega _{2}
\right\}\\
 \hphantom{\overline{C_{1;0}^{\left( 0\right) }} }{} = \big\{ \left( 0,0;1;0\right) ,\left( 2/3,1/3;1;0\right) ,\left(
1/3,2/3;1;0\right) \big\} ,
\end{gather*}
$\omega _{i}$ are the fundamental weights.

The classical components $\overset{\circ }{f\psi }$ of the folded fan shifts
\[
w_{\phi \left( \xi ,w\right) }\circ \left( \xi -\left( w\circ \rho -\rho
\right) \right) -\xi ,\qquad \xi \in \overline{C_{k}^{\left( 0\right) }}
\]
belong to the classical root lattice $Q\big( \overset{\circ }{\frak{g}}%
\big) $. For any weight $\xi =\big( \overset{\circ }{\xi };1;0\big) \in
\overline{C_{1;0}^{\left( 0\right) }}$ these classical components are equal
to zero, thus the folded fan has the form
\[
F\Psi \big( \overset{\circ }{\xi _{j}}\big) :=\big\{ \big(
0;0;n_{f\psi  ( \overset{\circ }{\xi } ) }\big) \big\} ,\qquad
\xi_{j} \in \overline{C_{k}^{\left( 0\right) }},\qquad j=1,2,3.
\]
It is convenient to indicate the multiplicities
\[
\eta _{j,s}\left( n\right) =
-\sum_{\substack{ \tilde{w}_{j,s}\in W,  \\
w_{\phi\left( \xi ,w\right) }\circ
\left( \xi _{j}-\left( \tilde{w}_{j,s}\circ \rho -\rho \right)
\right)
= ( \overset{\circ }{\xi _{s}};k;n_{j}+n )  }}
\epsilon (\tilde{w}_{j,s})
\]
as the additional coordinates of the shifts $f\psi $:
\[
F\Psi\big( \overset{\circ }{\xi _{j}}\big) :=\big\{ \big(
0,0,n_{f\psi  ( \overset{\circ }{\xi } ) },\eta _{j,s}\big(
n_{f\psi  ( \overset{\circ }{\xi } ) }\big) \big) \big\} .
\]
Thus any folded fan for the highest weight $\mu $ of level $k=1$ contains
only ``one string''. Moreover the fans $F\Psi\big( \overset{%
\circ }{\xi _{j}}\big) $ do not depend on the choice of $\xi _{j}=\big(
\overset{\circ }{\xi _{j}};1;0\big) \in \overline{C_{1;0}^{\left( 0\right)
}}$. The latter results are in full accord with the Proposition 12.6
in \cite{Kac}.

Using the fan $\Gamma _{\frak{h}\subset \frak{g}}$ we obtain the folded fan
(only the shifts with nonzero multiplicities $\eta _{j,j}$ are indicated,
the maximal grade here is $n=20$):
\begin{gather*}
F\Psi\big( \overset{\circ }{\xi _{j}}\big) :=\big\{
  \left( 0;0;0;-1\right) ,\left( 0;0;1;2\right) ,\left(
0;0;2;1\right) ,\left( 0;0;3;-2\right) ,\left( 0;0;4;-1\right) , \\
\hphantom{F\Psi\big( \overset{\circ }{\xi _{j}}\big) :=\big\{}{} \left(0;0;5;-2\right) ,\left( 0;0;7;2\right) ,\left( 0;0;8;2\right) ,
\left( 0;0;9;-1\right)
,\left( 0;0;10;-1\right) ,\left( 0;0;13;-2\right) , \\
\hphantom{F\Psi\big( \overset{\circ }{\xi _{j}}\big) :=\big\{}{} \left( 0;0;14;-3\right) ,\left( 0;0;15;2\right) ,\left( 0;0;16;-2\right) ,
\left( 0;0;19;2\right)
,\left( 0;0;20;2\right) ,\ldots  \big\}
\end{gather*}
The multiplicities
\begin{gather*}
\left\{ \eta _{j,j}\left( n\right) \right\} _{n=0,\ldots ,20}
 = \big\{ -1,2,1,-2,-1,-2,2,0,2,2,-1,-1,0,0,
  -2,-3,2,-2,0,0,2,2\big\}
\end{gather*}
form the unique nonzero matrix $\mathbf{M}_{\left( j,j\right) }$ for $j=1,2,3 $:
\[
\mathbf{M}_{\left( j,j\right) }:=
\begin{array}{cccc}
\eta _{j,j}\left( 0\right) & \eta _{j,j}\left( 1\right) & \cdots &
\eta _{j,j}\left( -n_{j}\right) \\
0 & \eta _{j,j}\left( 0\right) & \cdots & \eta _{j,j}\left( -n_{j} -1\right) \\
\vdots & \vdots & \vdots & \vdots \\
0 & 0 & \cdots & \eta _{j,j}\left( 0\right)
\end{array}
.
\]
The matrix $\mathbf{M}$ is block-diagonal and the equation (\ref{finite-string-solution}) splits into three equivalent (for $\mu =\left(
0,0;1;0\right) ,\left( 2/3,1/3;1;0\right) ,\left( 1/3,2/3;1;0\right) $)
relations $\mathbf{m}_{\left( j;-20\right) }^{\left( \mu \right) }=\mathbf{M}%
_{\left( j,j\right) }^{-1} \mathbf{\delta }_{\left( j;-20\right) }^{\mu }$
determining the unique string function with coef\/f\/icients $\mathbf{m}_{\left(
j;-20\right) }^{\left( \mu \right) }=\big( m_{j,-20}^{\left( \mu \right)
},m_{j,-19}^{\left( \mu \right) },\ldots ,m_{j,0}^{\left( \mu \right)
}\big) $,
\begin{gather*}
\sigma \left( q\right)  = 1+2q+5q^{2}+10q^{3}+20q^{4}+36q^{5}
 +65q^{6}+110q^{7}+185q^{8}+300q^{9}+481q^{10} \\
\hphantom{\sigma \left( q\right)  =}{} +752q^{11}+1165q^{12}+1770q^{13}+2665q^{14}+3956q^{15}
 +5822q^{16}+8470q^{17}\\
 \hphantom{\sigma \left( q\right)  =}{}
 +12230q^{18}+17490q^{19}+24842q^{20}+\cdots .
\end{gather*}
The obtained expression coincides with the expansion of the square of the
inverse Euler function (see Proposition~12.13 in~\cite{Kac} and the
relation~(12.13.4) there).

\subsubsection[$k=2$]{$\boldsymbol{k=2}$}\label{section5.1.2}

The set $\overline{C_{1;0}^{\left( 0\right) }}$ contains six weights:
\begin{gather*}
\overline{C_{1;0}^{\left( 0\right) }}  = \left\{
\begin{array}{c}
\left( 0,0;2;0\right) , ( \overset{\circ }{\omega _{1}};2;0 )
, ( \overset{\circ }{\omega _{2}};2;0 ) , \\
 ( \overset{\circ }{\omega _{1}}+\overset{\circ }{\omega _{2}}%
;2;0 ) , ( 2\overset{\circ }{\omega _{1}};2;0 ) , ( 2%
\overset{\circ }{\omega _{2}};2;0 )
\end{array}
\right\} = \\
\hphantom{\overline{C_{1;0}^{\left( 0\right) }} }{}  = \left\{
\begin{array}{c}
\left( 0,0;2;0\right) ,\left( 2/3,1/3;2;0\right) ,\left( 1/3,2/3;2;0\right) ,
 \\
 \left( 1,1;2;0\right) ,\left( 4/3,2/3;2;0\right) ,\left( 2/3,4/3;2;0\right)
\end{array}
\right\} .
\end{gather*}
This set is divided into 3 congruence classes. The fan shifts cannot connect
vectors from dif\/ferent classes. Thus instead of the set $\Xi _{2}$ we
can consider three subsets separately:
\begin{gather*}
\Xi _{\rm 2;I} =\left\{ \left( 0,0;2;0\right) ,\left( 1,1;2;0\right) \right\} ,
\\
\Xi _{\rm 2;II} =\left\{ \left( 2/3,1/3;2;0\right) ,\left( 2/3,4/3;2;0\right)
\right\} , \\
\Xi _{\rm 2;III} =\left\{ \left( 1/3,2/3;2;0\right) ,\left( 4/3,2/3;2;0\right)
\right\} .
\end{gather*}

Let us start with $\overset{\circ }{\xi _{s}}\in \Xi _{\rm 2;I}$ and $\mu
=\left( 0,0;2;0\right) $. Here we have two folded fans $F\Psi%
\big( \overset{\circ }{\xi _{1}}\big) $ and $F\Psi\big(
\overset{\circ }{\xi _{2}}\big) $.

Using the fan $\Gamma _{\frak{h}\subset \frak{g}}$ (\ref{star-a2-1}) we
obtain the folded fans (the maximal grade here is $n=9$):
\begin{gather*}
F\Psi\big( \overset{\circ }{\xi _{1}}\big) :=\big\{
  \left( 0;0;0;-1\right), \left( 0;0;2;1\right) ,\left( 0;0;4;2\right) ,\left(
0;0;8;-2\right) ,\left( 0;0;10;-2\right) , \\
\hphantom{F\Psi\big( \overset{\circ }{\xi _{1}}\big) :=\big\{}{} \left( 1;1;0;2\right) ,\left( 1;1;1;-1\right) ,\left( 1;1;2;-2\right)
,\left( 1;1;3;-2\right) ,\left( 1;1;4;2\right) , \\
\hphantom{F\Psi\big( \overset{\circ }{\xi _{1}}\big) :=\big\{}{}\left( 1;1;5;1\right) ,\left( 1;1;6;-2\right) ,\left( 1;1;7;2\right) ,\left(
1;1;9;-1\right) ,\ldots  \big\},
\\
F\Psi\big( \overset{\circ }{\xi _{2}}\big) :=\big\{
 \left( 0;0;1;1\right) ,\left( 0;0;3;-2\right) ,\left(
0;0;7;1\right) ,\left( 0;0;9;-2\right) , \\
\hphantom{F\Psi\big( \overset{\circ }{\xi _{2}}\big) :=\big\{}{} \left( 1;1;1;2\right) ,\left( 1;1;2;-2\right) ,\left( 1;1;4;1\right) ,\left(
1;1;5;2\right) ,\left( 1;1;6;2\right) , \\
\hphantom{F\Psi\big( \overset{\circ }{\xi _{2}}\big) :=\big\{}{}\left( 1;1;7;-2\right) ,\left( 1;1;8;-2\right) ,\left( 1;1;9;-2\right)
,\ldots  \big\}.
\end{gather*}
The multiplicities ($n=0,\ldots ,10$)
\begin{gather*}
\left\{ \eta _{1,1}\left( -10+n\right) \right\}   = \left\{
-1,0,1,0,2,0,0,0,-2,0,-2\right\} , \\
\left\{ \eta _{1,2}\left( -10+n\right) \right\}   = \left\{
2,-1,-2,-2,2,1,-2,2,0,-1,0\right\} , \\
\left\{ \eta _{2,1}\left( -10+n\right) \right\}   = \left\{
0,1,0,-2,0,0,0,1,0,-1,0\right\} , \\
\left\{ \eta _{2,2}\left( -10+n\right) \right\}   = \left\{
-1,2,-2,0,1,2,2,-2,-2,-2,0\right\} ,
\end{gather*}
form the matrices $\mathbf{M}_{\left( s,t\right) }^{\left( \Xi ,2;v\right) }$
for $s,t=1,2$:
\[
\mathbf{M}_{\left( s,t\right) }^{\left( \Xi ,2;v\right) }:=
\begin{array}{cccc}
\eta _{s,t}\left( 0\right)  & \eta _{s,t}\left( 1\right)  & \cdots  &
\eta _{s,t}\left( 10\right)  \\
0 & \eta _{s,t}\left( 0\right)  & \cdots  & \eta _{s,t}\left( 9\right)
\\
\vdots  & \vdots  & \vdots  & \vdots  \\
0 & 0 & \cdots  & \eta _{s,t}\left( 0\right)
\end{array}
.
\]
The block-matrix $\mathbf{M}^{\left( \Xi ,2;v\right) }$ is
\[
\mathbf{M}^{\left( \Xi ,2;v\right) }:=\left\|
\begin{array}{cc}
\mathbf{M}_{\left( 1,1\right) }^{\left( \Xi ,2;v\right) } & \mathbf{M}%
_{\left( 1,2\right) }^{\left( \Xi ,2;v\right) } \\
\mathbf{M}_{\left( 2,1\right) }^{\left( \Xi ,2;v\right) } & \mathbf{M}%
_{\left( 2,2\right) }^{\left( \Xi ,2;v\right) }
\end{array}
\right\| .
\]
The equation
\[
\mathbf{m}_{\left( -10\right) }^{\left( 0,0;2;0\right) }=\big( \mathbf{M}%
^{\left( \Xi ,2;v\right) }\big) ^{-1} \mathbf{\delta }_{\left(
-10\right) }^{\left( 0,0;2;0\right) }
\]
gives two string functions $\sigma _{\left( s;-10\right) }^{\left(
0,0;2;0\right) }$ with the coef\/f\/icients in the subsections of the vector $%
\mathbf{m}_{\left( -10\right) }^{\left( 0,0;2;0\right) }$:
\begin{gather*}
\sigma _{\left( 1;-10\right) }^{\left( 0,0;2;0\right) }
 = 1+2q+8q^{2}+20q^{3}+52q^{4}+116q^{5}  \\
 \hphantom{\sigma _{\left( 1;-10\right) }^{\left( 0,0;2;0\right) }=}{} +256q^{6}+522q^{7}+1045q^{8}+1996q^{9}+3736q^{10}+\cdots ,
\\
\sigma _{\left( 2;-10\right) }^{\left( 0,0;2;0\right) } =q+4q^{2}+12q^{3}+32q^{4}+77q^{5}  \\
\hphantom{\sigma _{\left( 2;-10\right) }^{\left( 0,0;2;0\right) } =}{} +172q^{6}+365q^{7}+740q^{8}+1445q^{9}+2736q^{10}+\cdots .
\end{gather*}

In the second congruence class $\Xi _{\rm 2;{\rm II}}{=}\left\{ \left(
2/3,1/3;2;0\right) ,\left( 2/3,4/3;2;0\right) \right\} $ put $\mu {=}\left(
2/3,1/3;2;0\right) $. Again we have two folded fans $F\Psi%
\big( \overset{\circ }{\xi _{1}}\big) $ and $F\Psi\big(
\overset{\circ }{\xi _{2}}\big) $.

The multiplicities ($n=0,\ldots ,10$):
\begin{gather*}
\left\{ \eta _{1,1}\left( -10+n\right) \right\}   = \left\{
-1,2,-2,0,1,2,2,-2,-2,-2,0\right\} , \\
\left\{ \eta _{1,2}\left( -10+n\right) \right\}   = \left\{
1,0,-2,0,0,0,1,0,-1,0,2\right\} , \\
\left\{ \eta _{2,1}\left( -10+n\right) \right\}   = \left\{
0,2,-1,-2,-2,2,1,-2,2,0,-1\right\} , \\
\left\{ \eta _{2,2}\left( -10+n\right) \right\}   = \left\{
-1,0,1,0,2,0,0,0,-2,0,-2\right\} .
\end{gather*}
form the matrices $\mathbf{M}_{\left( s,t\right) }^{\Xi 2;{\rm II}}$ for $s,t=1,2$
and the $22\times 22$ block-matrix $\mathbf{M}$
\[
\mathbf{M}^{\Xi 2;{\rm II}}:=\left\|
\begin{array}{cc}
\mathbf{M}_{\left( 1,1\right) }^{\Xi 2;{\rm II}} & \mathbf{M}_{\left( 1,2\right)
}^{\Xi 2;{\rm II}} \\
\mathbf{M}_{\left( 2,1\right) }^{\Xi 2;{\rm II}} & \mathbf{M}_{\left( 2,2\right)
}^{\Xi 2;{\rm II}}
\end{array}
\right\| .
\]
The equation
\[
\mathbf{m}_{\left( -10\right) }^{\left( 2/3,1/3;2;0\right) }=\left( \mathbf{M%
}^{\Xi 2;{\rm II}}\right) ^{-1} \mathbf{\delta }_{\left( -10\right) }^{\left(
2/3,1/3;2;0\right) }
\]
gives two string functions $\sigma _{\left( s;-10\right) }^{\left(
2/3,1/3;2;0\right) }$ for the module $L^{\left( 2/3,1/3;2;0\right) }$ with
the coef\/f\/icients in the subsections of the vector $\mathbf{m}_{\left(
-10\right) }^{\left( 2/3,1/3;2;0\right) }$:
\begin{gather*}
\sigma _{\left( 1;-10\right) }^{\left( 2/3,1/3;2;0\right) }
=1+4q+13q^{2}+36q^{3}+89q^{4}+204q^{5}  \\
\hphantom{\sigma _{\left( 1;-10\right) }^{\left( 2/3,1/3;2;0\right) }=}{} +441q^{6}+908q^{7}+1798q^{8}+3444q^{9}+6410q^{10}+\cdots ,
\\
\sigma _{\left( 2;-10\right) }^{\left( 2/3,1/3;2;0\right) }
= 2q+7q^{2}+22q^{3}+56q^{4}+136q^{5}  \\
\hphantom{\sigma _{\left( 2;-10\right) }^{\left( 2/3,1/3;2;0\right) }=}{} +300q^{6}+636q^{7}+1280q^{8}+2498q^{9}+4708q^{10}+\cdots .
\end{gather*}

For the third congruence class $\Xi _{2;{\rm III}}{=}\left\{ \left(
1/3,2/3;2;0\right) ,\left( 4/3,2/3;2;0\right) \right\} $ the folded fans $%
F\Psi\big( \overset{\circ }{\xi _{1}}\big) $ and
$F\Psi \big( \overset{\circ }{\xi _{2}}\big) $ are the same as for the
second one. As a result the string functions also coincide: $\sigma _{\left(
s;-10\right) }^{\left( 1/3,2/3;2;0\right) }=\sigma _{\left( s;-10\right)
}^{\left( 2/3,1/3;2;0\right) } $ in accord with the $A_{2}$ external
automorphism.

\subsubsection[$k=4$]{$\boldsymbol{k=4}$}\label{section5.1.3}

The set $\overline{C_{1;0}^{\left( 0\right) }}$ contains 15 projected
maximal weights
\begin{gather*}
\left\{ \xi _{j}\mid \xi _{j}\in \Xi_{4}
;j=1,\ldots ,p_{\max }=15\right\} ,
\\
\overline{C_{1;0}^{\left( 0\right) }}=\left\{
\begin{array}{c}
4\omega _{0},3\omega _{0}+\omega _{1},3\omega _{0}+\omega _{2},2\omega
_{0}+2\omega _{1},2\omega _{0}+2\omega _{2}, \\
2\omega _{0}+\omega _{1}+\omega _{2},\omega _{0}+3\omega _{1},\omega
_{0}+3\omega _{2},\omega _{0}+2\omega _{1}+\omega _{2}, \\
\omega _{0}+\omega _{1}+2\omega _{2},3\omega _{1}+\omega _{2},\omega
_{1}+3\omega _{2},2\omega _{1}+2\omega _{2},4\omega _{1},4\omega _{2}
\end{array}
\right\} .
\end{gather*}
This set is divided into 3 congruence classes. Instead of the set $\Xi
_{4}$ we can consider separately three subsets:
\begin{gather*}
\Xi _{\rm 4;I}  = \left\{ \left( 0,0;4;0\right) ,\left( 1,1;4;0\right) ,\left(
1,2;4;0\right) ,\left( 2,1;4;0\right) ,\left( 2,2;4;0\right)  \right\} , \\
\Xi _{\rm 4;II}  = \big\{
\left( 2/3,1/3;4;0\right) ,\left( 2/3,4/3;4;0\right) ,\left(
5/3,4/3;4;0\right) ,
\left( 5/3,7/3;4;0\right) ,\left( 8/3,4/3;4;0\right)
\big\} , \\
\Xi _{\rm 4;III}  = \big\{
\left( 1/3,2/3;4;0\right) ,\left( 4/3,2/3;4;0\right) ,\left(
4/3,5/3;4;0\right) ,
\left( 7/3,5/3;4;0\right) ,\left( 4/3,8/3;4;0\right)
\big\} .
\end{gather*}

Let us start with $\overset{\circ }{\xi _{s}}\in \Xi _{\rm 4;I}$ and $\mu
=\left( 0,0;4;0\right) $. Here we have 5 folded fans $F\Psi%
\big( \overset{\circ }{\xi _{s}}\big) $, $s=1,\ldots ,5$.

Using the fan $\Gamma _{\frak{h}\subset \frak{g}}$ (\ref{star-a2-1}) we
construct the folded fans (the maximal grade here is chosen to be $n=9$):
\begin{gather*}
F\Psi\big( \overset{\circ }{\xi _{1}}\big) :=\big\{
  \left( 0,0;0;-1\right) ,\left( 0,0;9;2\right)
,(1,1;0;2),(1,1;1;1),(1,1;3;-1), \\
\hphantom{F\Psi\big( \overset{\circ }{\xi _{1}}\big) :=\big\{}{}
(1,1;4;-2),(1,1;5;2),(1,1;6;-2),(1,1;7;-1),(1,1;8;2), \\
\hphantom{F\Psi\big( \overset{\circ }{\xi _{1}}\big) :=\big\{}{}(1,2;0;-1),(1,2;1;-1),(1,2;3;1),(1,2;5;1),(1,2;8;1), \\
\hphantom{F\Psi\big( \overset{\circ }{\xi _{1}}\big) :=\big\{}{}
(2,1;0;-1),(2,1;1;-1),(2,1;3;1),(2,1;5;1),(2,1;8;1), \\
\hphantom{F\Psi\big( \overset{\circ }{\xi _{1}}\big) :=\big\{}{}
(2,2;0;1),(2,2;2;2),(2,2;4;-2),(2,2;6;-2),(2,2;8;-2),
\ldots \big\},
\\
F\Psi\big( \overset{\circ }{\xi _{2}}\big) :=\big\{
 \left( 0,0;1;1\right) ,\left( 0;0;5;-1\right) ,\left(
1;1;0;-1\right) ,\left( 1;1;2;-1\right) ,\left( 1;1;4;2\right) , \left( 1;1;5;-2\right) ,\\
\hphantom{F\Psi\big( \overset{\circ }{\xi _{2}}\big) :=\big\{}{}
\left( 1;1;8;2\right) ,\left( 1;1;9;2\right) ,\left(
1;1;1;2\right) ,\left( 1;1;2;-2\right) , \left( 1;1;4;1\right) ,\left( 1;1;5;2\right) ,\\
\hphantom{F\Psi\big( \overset{\circ }{\xi _{2}}\big) :=\big\{}{}
\left( 1;1;6;2\right) ,\left(
1;2;0;1\right) ,\left( 1;2;1;-1\right) ,\left( 1;2;2;1\right) ,\left( 1;2;4;1\right) ,\left( 1;2;5;-1\right) , \\
\hphantom{F\Psi\big( \overset{\circ }{\xi _{2}}\big) :=\big\{}{}
\left(1;2;6;-1\right) ,\left( 1;2;7;-1\right) , \left( 2;1;0;1\right) ,\left( 2;1;1;-1\right) ,\left( 2;1;2;1\right) ,\left(2;1;4;1\right) ,\\
\hphantom{F\Psi\big( \overset{\circ }{\xi _{2}}\big) :=\big\{}{}
\left( 2;1;5;-1\right) ,
\left( 2;1;6;-1\right) ,\left( 2;1;7;-1\right) ,\left( 2;2;0;1\right)
,\left( 2;2;2;2\right) ,\left( 2;2;4;-2\right) , \\
\hphantom{F\Psi\big( \overset{\circ }{\xi _{2}}\big) :=\big\{}{}
\left( 2;2;6;-2\right) ,\left( 2;2;8;-2\right) ,\ldots \big\},
\\
F\Psi\big( \overset{\circ }{\xi _{3}}\big) :=\big\{
  \left( 0,0;2;-1\right) ,\left( 0;0;6;1\right) ,\left(
1;1;1;1\right) ,\left( 1;1;4;2\right) , \left( 1;1;6;-2\right) ,\left( 1;1;7;-2\right) ,\\
\hphantom{F\Psi\big( \overset{\circ }{\xi _{3}}\big) :=\big\{}{}
\left( 1;2;0;-1\right)
,\left( 1;2;1;1\right) ,\left( 1;2;4;-1\right) ,\left( 1;1;5;-1\right) ,\left( 1;2;6;1\right) ,\left( 1;2;8;1\right) ,\\
\hphantom{F\Psi\big( \overset{\circ }{\xi _{3}}\big) :=\big\{}{}
\left(
1;2;9;2\right) , \left( 2;1;1;-1\right) , \left( 2;1;2;-1\right) ,\left( 2;1;3;2\right) ,\left( 2;1;5;1\right) ,
\left(2;1;6;-1\right) ,\\
\hphantom{F\Psi\big( \overset{\circ }{\xi _{3}}\big) :=\big\{}{}
\left( 2;1;8;-1\right) ,
\left( 2;2;0;1\right) ,\left( 2;2;2;-1\right) ,\left( 2;2;8;1\right) ,\ldots
\big\}.
\end{gather*}
The fan $F\Psi\big( \overset{\circ }{\xi _{4}}\big) $ is
equal to $\big\{ F\Psi\big( \overset{\circ }{\xi _{3}}\big)
|\left( 1;2;n;m\right) \rightleftharpoons \left( 2;1;n;m\right) \big\} $
\begin{gather*}
F\Psi\big( \overset{\circ }{\xi _{5}}\big) :=\big\{
 \left( 0,0;4;1\right) ,\left( 0;0;8;-2\right) ,\left(
1;1;1;1\right) ,\left( 1;1;2;-2\right) , \\
\hphantom{F\Psi\big( \overset{\circ }{\xi _{5}}\big) :=\big\{}{}
\left( 1;1;3;-2\right) ,\left( 1;1;4;2\right) ,\left( 1;1;5;1\right) ,\left(
1;1;7;-1\right) ,\left( 1;1;8;2\right) , \\
\hphantom{F\Psi\big( \overset{\circ }{\xi _{5}}\big) :=\big\{}{}
\left( 1;2;1;1\right) ,\left( 1;2;2;-1\right) ,\left( 1;2;6;-1\right)
,\left( 1;2;7;1\right) ,\left( 1;2;9;1\right) , \\
\hphantom{F\Psi\big( \overset{\circ }{\xi _{5}}\big) :=\big\{}{}
\left( 2;1;1;1\right) ,\left( 2;1;2;-1\right) ,\left( 2;1;6;-1\right)
,\left( 2;1;7;1\right) ,\left( 2;1;9;1\right) , \\
\hphantom{F\Psi\big( \overset{\circ }{\xi _{5}}\big) :=\big\{}{}\left( 2;2;0;-1\right) ,\left( 2;2;2;2\right) ,\left( 2;2;6;-2\right)
,\left( 2;2;8;-1\right) ,\ldots \big\}.
\end{gather*}

Their multiplicities (for $n=0,\ldots ,9$)
\begin{gather*}
\left\{ \eta _{1,1}\left( -9+n\right) \right\} =\left\{
-1,0,0,0,0,0,0,0,2,0\right\} , \\
\left\{ \eta _{1,2}\left( -9+n\right) \right\} =\left\{
2,1,0,-1,-2,2,-2,-1,2,0\right\} , \\
\left\{ \eta _{1,3}\left( -9+n\right) \right\} =\left\{
-1,-1,0,1,0,1,0,0,1,0\right\} , \\
\left\{ \eta _{1,4}\left( -9+n\right) \right\} =\left\{
-1,-1,0,1,0,1,0,0,1,0\right\} , \\
\left\{ \eta _{1,5}\left( -9+n\right) \right\} =\left\{
1,0,2,0,-2,0,-2,0,-2,0\right\} ,
\\
\left\{ \eta _{2,1}\left( -9+n\right) \right\} =\left\{
0,1,0,0,0,-1,0,0,0,0\right\} , \\
\left\{ \eta _{2,2}\left( -9+n\right) \right\} =\left\{
-1,0,-1,0,2,-2,0,0,2,2\right\} , \\
\left\{ \eta _{2,3}\left( -9+n\right) \right\} =\left\{
1,-1,1,0,1,-1,-1,-1,0,0\right\} , \\
\left\{ \eta _{2,4}\left( -9+n\right) \right\} =\left\{
1,-1,1,0,1,-1,-1,-1,0,0\right\} , \\
\left\{ \eta _{2,5}\left( -9+n\right) \right\} =\left\{
1,0,2,0,-2,0,-2,0,-2,0\right\} ,
\\
\left\{ \eta _{3,1}\left( -9+n\right) \right\}=\left\{
0,0,-1,0,0,0,1,0,0,0\right\} , \\
\left\{ \eta _{3,2}\left( -9+n\right) \right\} =\left\{
0,1,0,0,2,0,-2,-2,0,0\right\} , \\
\left\{ \eta _{3,3}\left( -9+n\right) \right\} =\left\{
-1,1,0,0,-1,-1,1,0,1,2\right\} , \\
\left\{ \eta _{3,4}\left( -9+n\right) \right\} =\left\{
0,-1,-1,2,0,1,-1,0,-1,0\right\} , \\
\left\{ \eta _{3,5}\left( -9+n\right) \right\} =\left\{
1,0,-1,0,0,0,0,0,1,0\right\} ,
\\
\left\{ \eta _{4,1}\left( -9+n\right) \right\} =\left\{
0,0,-1,0,0,0,1,0,0,0\right\} , \\
\left\{ \eta _{4,2}\left( -9+n\right) \right\} =\left\{
0,1,0,0,2,0,-2,-2,0,0\right\} , \\
\left\{ \eta _{4,3}\left( -9+n\right) \right\} =\left\{
0,-1,-1,2,0,1,-1,0,-1,0\right\} , \\
\left\{ \eta _{4,4}\left( -9+n\right) \right\} =\left\{
-1,1,0,0,-1,-1,1,0,1,2\right\} , \\
\left\{ \eta _{4,5}\left( -9+n\right) \right\} =\left\{
1,0,-1,0,0,0,0,0,1,0\right\} ,
\\
\left\{ \eta _{5,1}\left( -9+n\right) \right\} =\left\{
0,0,0,0,1,0,0,0,-2,0\right\} , \\
\left\{ \eta _{5,2}\left( -9+n\right) \right\} =\left\{
0,1,-2,-2,2,1,0,-1,2,0\right\} , \\
\left\{ \eta _{5,3}\left( -9+n\right) \right\} =\left\{
0,1,-1,0,0,0,-1,1,0,1\right\} , \\
\left\{ \eta _{5,4}\left( -9+n\right) \right\} =\left\{
0,1,-1,0,0,0,-1,1,0,1\right\} , \\
\left\{ \eta _{5,5}\left( -9+n\right) \right\} =\left\{
-1,0,2,0,0,0,-2,0,-1,0\right\} .
\end{gather*}

The matrices $\mathbf{M}_{\left( s,t\right) }^{\Xi 4;{\rm I}}$ for $s,t=1,\ldots ,5
$:
\[
\mathbf{M}_{\left( s,t\right) }^{\Xi 4;{\rm I}}:=
\begin{array}{cccc}
\eta _{s,t}\left( 0\right)  & \eta _{s,t}\left( 1\right)  & \cdots  & \eta
_{s,t}\left( 9\right)  \\
0 & \eta _{s,t}\left( 0\right)  & \cdots  & \eta _{s,t}\left( 8\right)  \\
\vdots  & \vdots  & \vdots  & \vdots  \\
0 & 0 & \cdots  & \eta _{s,t}\left( 0\right)
\end{array}
.
\]
For example,
\[
\mathbf{M}_{\left( 2,3\right) }^{\Xi 4;{\rm I}}:=
\begin{array}{rrrrrrrrrr}
1 & -1 & 1 & 0 & 1 & -1 & -1 & -1 & 0 & 0 \\
0 & 1 & -1 & 1 & 0 & 1 & -1 & -1 & -1 & 0 \\
0 & 0 & 1 & -1 & 1 & 0 & 1 & -1 & -1 & -1 \\
0 & 0 & 0 & 1 & -1 & 1 & 0 & 1 & -1 & -1 \\
0 & 0 & 0 & 0 & 1 & -1 & 1 & 0 & 1 & -1 \\
0 & 0 & 0 & 0 & 0 & 1 & -1 & 1 & 0 & 1 \\
0 & 0 & 0 & 0 & 0 & 0 & 1 & -1 & 1 & 0 \\
0 & 0 & 0 & 0 & 0 & 0 & 0 & 1 & -1 & 1 \\
0 & 0 & 0 & 0 & 0 & 0 & 0 & 0 & 1 & -1 \\
0 & 0 & 0 & 0 & 0 & 0 & 0 & 0 & 0 & 1
\end{array}
.
\]
Matrices $\mathbf{M}_{\left( s,t\right) }^{\Xi 4;{\rm I}}$ form the block-matrix $%
\mathbf{M}^{\Xi 4;{\rm I}}\mathbf{=}\left\| \mathbf{M}_{\left( s,t\right) }^{\Xi
4;{\rm I}}\right\| _{s,t=1,\ldots ,5}$. With this matrix we can describe f\/ive
modules of the level 4 with the highest weights
$\mu _{s}\in \Xi
_{4;{\rm I}}=\left\{ \left( 0,0;4;0\right) ,\left( 1,1;4;0\right),\right.$
$\left. \left(
1,2;4;0\right) ,\left( 2,1;4;0\right) ,\left( 2,2;4;0\right) \right\} $ . We
construct f\/ive sets of string functions $\sigma _{\left( t;-9\right)
}^{\left( \mu _{s}\right) }$ in terms of their coef\/f\/icients obtained as ten
dimensional subsections of the vector $\mathbf{m}_{\left( -9\right)
}^{\left( \mu _{s}\right) }$:
\[
\mathbf{m}_{\left( -9\right) }^{\left( \mu _{s}\right) }=\left( \mathbf{M}%
^{\Xi 4;{\rm I}}\right) ^{-1} \mathbf{\delta }_{\left( -9\right) }^{\left( \mu
_{s}\right) }.
\]
The answer is as follows:
\begin{gather*}
\sigma _{\left( 1;-9\right) }^{\left( 0,0;4;0\right) }
=1+2q+8q^{2}+24q^{3}+72q^{4}+190q^{5}+490q^{6}  \\
\hphantom{\sigma _{\left( 1;-9\right) }^{\left( 0,0;4;0\right) }=}{}
+1176q^{7}+2729q^{8}+6048q^{9}+\cdots , \\
\sigma _{\left( 2;-9\right) }^{\left( 0,0;4;0\right) }
 = q+4q^{2}+15q^{3}+48q^{4}+138q^{5}+366q^{6}
  +913q^{7}+2156q^{8}+4874q^{9}+\cdots , \\
\sigma _{\left( 3;-9\right) }^{\left( 0,0;4;0\right) }
 = q^{2}+6q^{3}+23q^{4}+74q^{5}+2121q^{6}+556q^{7}
  +1366q^{8}+3184q^{9}+\cdots , \\
\sigma _{\left( 4;-9\right) }^{\left( 0,0;4;0\right) }
 = q^{2}+6q^{3}+23q^{4}+74q^{5}+2121q^{6}+556q^{7}
  +1366q^{8}+3184q^{9}+\cdots , \\
\sigma _{\left( 5;-9\right) }^{\left( 0,0;4;0\right) }
 = q^{2}+4q^{3}+18q^{4}+56q^{5}+167q^{6}+440q^{7}
  +1103q^{8}+2588q^{9}+\cdots ,
\\
\sigma _{\left( 1;-9\right) }^{\left( 1,1;4;0\right) }
 = 2+10q+40q^{2}+133q^{3}+398q^{4}+1084q^{5}+2760q^{6} \\
\hphantom{\sigma _{\left( 1;-9\right) }^{\left( 1,1;4;0\right) }=}{} +6632q^{7}+15214q^{8}+33508q^{9}+\cdots , \\
\sigma _{\left( 2;-9\right) }^{\left( 1,1;4;0\right) }
 = 1+6q+27q^{2}+96q^{3}+298q^{4}+836q^{5}+2173q^{6} \\
 \hphantom{\sigma _{\left( 2;-9\right) }^{\left( 1,1;4;0\right) }=}{}  +5310q^{7}+12341q^{8}+27486q^{9}+\cdots , \\
\sigma _{\left( 3;-9\right) }^{\left( 1,1;4;0\right) }
 = 2q^{2}+12q^{3}+49q^{4}+166q^{5}+494q^{6}+1340q^{7}
  +3387q^{8}+8086q^{9}+\cdots , \\
\sigma _{\left( 4;-9\right) }^{\left( 1,1;4;0\right) }
 = 2q^{2}+12q^{3}+49q^{4}+166q^{5}+494q^{6}+1340q^{7}
  +3387q^{8}+8086q^{9}+\cdots , \\
\sigma _{\left( 5;-9\right) }^{\left( 1,1;4;0\right) }
 = q+8q^{2}+35q^{3}+124q^{4}+379q^{5}+1052q^{6}+2700q^{7}
  +6536q^{8}+15047q^{9}+\cdots ,
\\
\sigma _{\left( 1;-9\right) }^{\left( 1,2;4;0\right) }
 = 1+8q+32q^{2}+110q^{3}+322q^{4}+872q^{5}+2183q^{6}  \\
 \hphantom{\sigma _{\left( 1;-9\right) }^{\left( 1,2;4;0\right) }=}{}
  +5186q^{7}+11730q^{8}+25552q^{9}+\cdots , \\
\sigma _{\left( 2;-9\right) }^{\left( 1,2;4;0\right) }
 = 1+6q+25q^{2}+85q^{3}+255q^{4}+695q^{5}+1764q^{6} \\
\hphantom{\sigma _{\left( 2;-9\right) }^{\left( 1,2;4;0\right) }=}{} +4226q^{7}+9653q^{8}+21179q^{9}+\cdots , \\
\sigma _{\left( 3;-9\right) }^{\left( 1,2;4;0\right) }
 = 1+4q+16q^{2}+54q^{3}+163q^{4}+450q^{5}+1161q^{6}+2824q^{7} \\
\hphantom{\sigma _{\left( 3;-9\right) }^{\left( 1,2;4;0\right) }=}{}  +6549q^{8}+14572q^{9}+\cdots , \\
\sigma _{\left( 4;-9\right) }^{\left( 1,2;4;0\right) }
 = 2q+11q^{2}+44q^{3}+143q^{4}+414q^{5}+1096q^{6}+2714q^{7}  \\
 \hphantom{\sigma _{\left( 4;-9\right) }^{\left( 1,2;4;0\right) }=}{}
  +6364q^{8}+14272q^{9}+\cdots , \\
\sigma _{\left( 5;-9\right) }^{\left( 1,2;4;0\right) }
 = 2q+9q^{2}+36q^{3}+115q^{4}+336q^{5}+890q^{6}+2224q^{7}
  +5241q^{8}+11840q^{9}+\cdots .
\end{gather*}
The next set of string functions $\sigma _{\left( s;-9\right) }^{\left(
2,1;4;0\right) }$ coincides with the previous one where the third and the
fourth strings are interchanged: $\sigma _{\left( 3;-9\right) }^{\left(
2,1;4;0\right) }=\sigma _{\left( 4;-9\right) }^{\left( 1,2;4;0\right)
},\sigma _{\left( 4;-9\right) }^{\left( 2,1;4;0\right) }=\sigma _{\left(
3;-9\right) }^{\left( 1,2;4;0\right) }.$ The last set describes the module
$L^{\mu _{5}}$ where $\mu _{5}$ is the highest weight in $\Xi _{4;{\rm I}}$:
\begin{gather*}
\sigma _{\left( 1;-9\right) }^{\left( 2,2;4;0\right) }
 = 3+14q+58q^{2}+184q^{3}+536q^{4}+1408q^{5}+3492q^{6} \\
 \hphantom{\sigma _{\left( 1;-9\right) }^{\left( 2,2;4;0\right) }=}{} +8160q^{7}+18299q^{8}+39428q^{9}+\cdots , \\
\sigma _{\left( 2;-9\right) }^{\left( 2,2;4;0\right) }
 = 2+11q+44q^{2}+145q^{3}+424q^{4}+1133q^{5}+2830q^{6} \\
\hphantom{\sigma _{\left( 2;-9\right) }^{\left( 2,2;4;0\right) }=}{}  +6688q^{7}+15102q^{8}+32805q^{9}+\cdots , \\
\sigma _{\left( 3;-9\right) }^{\left( 2,2;4;0\right) }
 = 1+6q+25q^{2}+86q^{3}+260q^{4}+716q^{5}+1833q^{6}+4426q^{7} \\
\hphantom{\sigma _{\left( 3;-9\right) }^{\left( 2,2;4;0\right) }=}{}  +10183q^{8}+22488q^{9}+\cdots , \\
\sigma _{\left( 4;-9\right) }^{\left( 2,2;4;0\right) }
 = 1+6q+25q^{2}+86q^{3}+260q^{4}+716q^{5}+1833q^{6}+4426q^{7} \\
\hphantom{\sigma _{\left( 4;-9\right) }^{\left( 2,2;4;0\right) }=}{}  +10183q^{8}+22488q^{9}+\cdots , \\
\sigma _{\left( 5;-9\right) }^{\left( 2,2;4;0\right) }
 = 1+4q+19q^{2}+64q^{3}+202q^{4}+560q^{5}+1464q^{6}+3568q^{7} \\
\hphantom{\sigma _{\left( 5;-9\right) }^{\left( 2,2;4;0\right) }=}{}  +8315q^{8}+18512q^{9}+\cdots .
\end{gather*}

Notice that in the congruence class $\Xi _{4;{\rm I}}$ we have only 17
dif\/ferent string functions.

\section{Conclusions}\label{section6}

The folded fans $F\Psi\big( \overset{\circ }{\xi _{j}}\big) $
(for a f\/ixed level $k$ and the congruence class $\Xi _{k;v}$ of weights in $
\overline{C_{k}^{\left( 0\right) }}$) were constructed by transporting
to the fundamental Weyl chamber the standard set $\widehat{%
\Psi ^{\left( 0\right) }}$~-- the set of singular weights of module $L^{0}$
supplied with the anomalous multiplicities. We have found out that the
shifts $f\psi \big( \overset{\circ }{\xi }\big) \in F\Psi %
\big( \overset{\circ }{\xi }\big) $ (connecting $\xi _{j}\in \Xi _{k;v}$)
together with their multiplicities~$\eta _{j,s}$ describe the recursive
properties of the weights of modules $L^{\xi _{j}}$ with the highest weights~$\xi _{j}$. Thus the set $\big\{ F\Psi\big( \overset{\circ }{%
\xi _{j}}\big) |\,\xi _{j}\in \Xi _{k;v}\big\} $ describes the recursive
properties of the string functions $\big\{ \sigma _{j}^{\mu ,k}\,|\,\mu ,\xi
_{j}\in \Xi _{k;v}\big\} $. When for a f\/ixed module $L^{\mu }$ these
properties are simultaneously
considered for $\big\{ \sigma _{j}^{\mu ,k}\,|\,\mu ,\xi _{j}\in
\Xi _{k;v}\big\} $ they can be written in a form of the
equation $\mathbf{M}^{\Xi ,k;v} \mathbf{m}_{\left( u\right) }^{\left( \mu
\right) }=\mathbf{\delta }_{\left( u\right) }^{\mu }$. In this equation $%
\mathbf{M}^{\Xi ,k;v}$ is a matrix formed by the multiplicities $\eta _{j,s}$
of the fan shifts, $\mathbf{\delta }_{\left( u\right) }^{\mu }$ indicates
what weight in the set $ \Xi _{k;v}$ is chosen to be the
highest weight $\mu$ of the module
and $\mathbf{m}_{\left( u\right) }^{\left( \mu \right) }$ is a
vector of string functions coef\/f\/icients. As far as $\mathbf{M}^{\Xi ,k;v}$
is invertible the solution $\mathbf{m}_{\left( u\right) }^{\left( \mu
\right) }=\left( \mathbf{M}^{\Xi ,k;v}\right) ^{-1} \mathbf{\delta }%
_{\left( u\right) }^{\mu }$ can be explicitly written and the full set of
string functions $\big\{ \sigma _{j}^{\mu ,k}\,|\, \mu ,\xi _{j}\in \Xi
_{k;v}\big\} $ for $L^{\mu }$ is determined by this linear equation (at
least for any common f\/inite ``length'' of all the strings).

There are two points that we want to stress.
The f\/irst is that in this algorithm the singular vectors
$\psi \in \widehat{\Psi ^{\left( \mu \right) }}$ of $L^{\mu }$ are not needed
(except the highest weight $\mu $). The second point is that the
crossections $F\Psi\big( \overset{\circ }{\xi _{j}}\big)
\cap \overline{C_{k,0}^{\left( 0\right) }}$ form the parts of the classical
folded fans for $\overset{\circ }{\frak{g}}$. It can be easily verif\/ied
that the string starting vectors $\big\{ \sigma _{j}^{\mu ,k}\,|\, \mu
,\xi _{j}\in \Xi _{k;v};\  n=0\big\} $ and their multiplicities present
the diagram $\mathcal{N}^{\overset{\circ }{\mu }}\cap \overline{C_{k}^{\left(
0\right) }}$  of the module $L^{\overset{\circ }{\mu }}\big( \overset{\circ
}{\frak{g}}\big) $. In general the crossections $F\Psi\big(
\overset{\circ }{\xi _{j}}\big) \cap \overline{C^{\left( 0\right) }\big(
\overset{\circ }{\frak{g}}\big) }$ do not coincide with the classical
folded fans because the chambers $\overline{C^{\left( 0\right) }\big(
\overset{\circ }{\frak{g}}\big) }$ are inf\/inite (contrary to~$\overline{%
C_{k,0}^{\left( 0\right) }}$ for any f\/inite $k$).

As it was demonstrated in the examples the folded fans provide an ef\/fective
tool when studying the string functions for integrable highest weight
modules of af\/f\/ine Lie algebras.

\subsection*{Acknowledgements}

The authors appreciate helpful remarks made by the Referees.
The work was supported in
part by RFBR grants N 06-01-00451, N 08-01-00638
and the National Project RNP.2.1.1.1112.

\newpage
\pdfbookmark[1]{References}{ref}
\LastPageEnding


\begin{thebibliography}{99}

\footnotesize\itemsep=0pt

\bibitem{DJKMO}   Date E.,  Jimbo M., Kuniba A., Miwa T., Okado M., One-dimensional conf\/iguration sums in vertex models and af\/f\/ine Lie algerba
characters,  {\it Lett. Math. Phys.} {\bf 17} (1989), 69--77 (Preprint RIMS-631, 1988).

\bibitem{LD}  Faddeev L.D., How the algebraic Bethe ansatz works for integrable models,
in  Symetries Quantique, Proc. Les Houches Summer School (Les Houches, 1995),
 Editors A.~Connes, K.~Gawedski, J.~Zinn-Justin, North-Holland,
1998, 149--219, \href{http://arxiv.org/abs/hep-th/9605187}{hep-th/9605187}.

\bibitem{KAZA}  Kazakov V.A., Zarembo K., Classical/quantum integrability in
non-compact sector of AdS/CFT,
{\it J. High Energy Phys.} {\bf 2004} (2004), no.~10, 060, 23 pages, \href{http://arxiv.org/abs/hep-th/0410105}{hep-th/0410105}.


\bibitem{BE}  Beisert N., The dilatation operator of $N=4$ super Yang--Mills  theory
and integrability, {\it Phys. Rep.} {\bf 405}  (2005), 1--202, \href{http://arxiv.org/abs/hep-th/0407277}{hep-th/0407277}.

\bibitem{BGG} Bernstein I.N., Gelfand I.M., Gelfand S.I., Dif\/ferential
operators on the basic af\/f\/ine space and a study of $\gamma $-modules, in Lie
Groups and Their Representations, Summer School of Bolyai Janos Math.
Soc. (Budapest, 1971), Editor I.M.~Gelfand, Halsted Press, New York,  1975, 21--64.

\bibitem{Kac}  Kac V., Inf\/inite-dimensional Lie algebras, 3rd ed.,
Cambridge University Press, Cambridge, 1990.

\bibitem{Wak1}  Wakimoto M., Inf\/inite-dimensional Lie algebras,
{\it Translations of Mathematical Monographs}, Vol.~195, American Mathematical Society, Providence, RI, 2001.

\bibitem{FauKing}  Fauser B., Jarvis P.D., King R.C., Wybourn B.G., New
branching rules induced by plethysm, \mbox{\href{http://arxiv.org/abs/math-ph/0505037}{math-ph/0505037}}.

\bibitem{Hwang}  Hwang S., Rhedin H., General branching functions of
af\/f\/ine Lie algebras, {\it Modern Phys. Lett.~A} {\bf 10} (1995), 823--830, \href{http://arxiv.org/abs/hep-th/9408087}{hep-th/9408087}.

\bibitem{FeigJimbo}  Feigin B., Feigin E., Jimbo M., Miwa T., Mukhin E.,
Principal $\widehat{\frak{sl_{3}}}$ subspaces and quantum Toda
Hamiltonians, \href{http://arxiv.org/abs/0707.1635}{arXiv:0707.1635}.

\bibitem{IKL-1}  Ilyin M., Kulish P., Lyakhovsky V., On a property of
branching coef\/f\/icients for af\/f\/ine Lie algebras, {\it Algebra i Analiz}, to appear, \href{http://arxiv.org/abs/0812.2124}{arXiv:0812.2124}.

\bibitem{Wak2}  Wakimoto M., Lectures on inf\/inite-dimensional Lie
algebra, World Scientif\/ic Publishing Co., Inc., River Edge, NJ, 2001.

\bibitem{Ful} Fulton W., Harris J., Representation theory. A f\/irst course,
{\it Graduate Texts in Mathematics}, Vol.~129, Springer-Verlag, New York, 1991.

\bibitem{LDu}  Lyakhovsky V.D., Recurrent properties of af\/f\/ine Lie algebra
representations, in Supersymmetry and Quantum Symmetry (August, 2007, Dubna),
to appear. 

\end{thebibliography}
\end{document}